\RequirePackage{fix-cm}
\documentclass[smallextended]{svjour3}       % onecolumn (second format)
\smartqed  % flush right qed marks, e.g. at end of proof
\usepackage{graphicx}
\usepackage[utf8]{inputenc}
\usepackage[english]{babel}
\usepackage{bbm}

\usepackage[nottoc]{tocbibind}
\usepackage{mathrsfs,amsfonts,amssymb,amsmath}

\usepackage{amsthm}
\usepackage{enumerate}
\usepackage{graphicx,cite}
\usepackage{romannum}
\usepackage{hyperref}
\usepackage{authblk}
\usepackage{graphicx}

\allowdisplaybreaks

\hypersetup{colorlinks=true, linkcolor=blue, citecolor=red}

\DeclareMathOperator{\diam}{diam}
\DeclareMathOperator{\Leb}{Leb}

\DeclareMathOperator{\interior}{int}

\DeclareMathOperator{\Emb}{Emb}

\newcommand{\Lip}{\mbox{\rm Lip}}

\numberwithin{equation}{section}

\newtheorem{assumption}{Assumption}

\usepackage[square,numbers]{natbib}
\begin{document}
\pagenumbering{arabic}
\setlength{\belowdisplayskip}{0pt}

\title{Statistics of Weakly Chaotic Systems}


\author{Leonid A. Bunimovich  \and
        Yaofeng Su
}


\institute{Leonid A. Bunimovich \at
              School of Mathematics, Georgia Institute of Technology,  Atlanta, USA\\
             \email{leonid.bunimovich@math.gatech.edu}           %  \\
           \and
            Yaofeng Su  \at
              Southern University of Science and Technology, Shenzhen, China\\
              \email{suyaofengmath@gmail.com} 
}

\date{Received: date / Accepted: date}

\maketitle

\begin{abstract}
One of the major breakthroughs in science of the last (20th) century was building a bridge between the worlds of
stochastic (random) systems and deterministic (dynamical) systems. It was started by the celebrated
1958 paper by A.N.Kolmogorov \cite{Kolmo}, who called this new theory (and actually a new way of thinking about deterministic systems) stochasticity of dynamical systems. Later, this name was essentially replaced by a short (sexier but more vague ``chaos theory"). Kolmogorov's discovery demonstrated that the time evolution of deterministic systems could be indistinguishable from the evolution of purely random (stochastic) systems. Moreover, it has been later established that typical deterministic systems are chaotic. However, as well as stochastic systems, which could be more or less random (from random processes with independent values to random processes with long memory) , chaotic dynamical systems could be strongly or weakly chaotic. Actually, the majority of chaotic systems are weakly chaotic, especially those that are relevant models for various real-world processes.   Naturally, the theory of weakly chaotic systems is less developed because it deals with more complicated problems than the studies of strongly chaotic systems. We present here a brief review of this theory, including recently published and some new results.

\end{abstract}

\tableofcontents

\section{Introduction}\ \par
In this paper we concentrate on weakly chaotic deterministic (dynamical) systems.
Therefore, everywhere we just verbally (without exact formulations) describe the results for
strongly chaotic systems and present the best (so far) results for weakly chaotic systems, which are more difficult to analyze because they occupy an intermediate position between strongly chaotic systems and dynamical systems with a non-chaotic, regular behavior. So, neither the methods for studying regular dynamics nor the ones for dealing with strongly chaotic systems are applicable to weakly chaotic systems.
We will most often consider examples of billiards and statistical physics systems to demonstrate what the weakly chaotic dynamical systems really are and where they naturally appear as relevant models of various processes in science (mostly in physics). In this review paper, we mean by weak chaos slow mixing rates, that is, if an orbit of the dynamics $f$ can densely fill any subregion of the phase space $\mathcal{M}$ but after a long time. More precisely, let $\phi, \psi$ be functions (observables) defined in $\mathcal{M}$ and let $\mu_\mathcal{M}$ be an invariant measure of dynamics (i.e. $\mu_{\mathcal{M}} \circ f^{-1}=\mu_{\mathcal{M}}$), then 

\[\lim_{n \to \infty}|\int \phi \circ f^n \psi d\mu_{\mathcal{M}}-\int \psi d\mu_{\mathcal{M}} \int \phi d\mu_{\mathcal{M}}|= 0.\]

The organization of the paper is the following. Section 2 deals with the strong law of large numbers (ergodic theorems) and large deviations (with emphasis on maximal-type large deviations, which indicates a core of principal differences between strongly and weakly chaotic systems). In Section 3 we consider the processes of recurrences and related Poisson approximations. In Section 4, we consider central limit theorems and the invariance principles. The stable limit laws are presented in Section 5.\\

We start by introducing some notations throughout the whole paper.
\begin{enumerate}
 \item $C_z$ denotes a constant depending on $z$.
    \item The notation $``a_n \precsim_z b_n"$  ($``a_n=O_{z}(b_n)"$) means that there is a constant $C_z \ge 1$ such that (s.t.) $ a_n \le C_z  b_n$ for all $n \ge 1$, whereas the notation $``a_n \precsim b_n"$ (or $``a_n=O(b_n)"$) means that there is a constant $C \ge 1$ such that $ a_n \le C  b_n$ for all $n \ge 1$. Next, $``a_n \approx_z b_n"$ and $a_n=C_z^{\pm 1}b_n$ mean that there is a constant $C_z \ge 1$ such that  $ C_z^{-1}  b_n \le a_n \le C_z b_n$ for all $n \ge 1$. Further, the notations $``a_n=C^{\pm1} b_n"$ and $``a_n \approx b_n"$ mean that there is a constant $C \ge 1$ such that $ C^{-1}  b_n \le a_n \le C b_n$ for all $n \ge 1$. Finally, $``a_n =o(b_n)"$ means that  $\lim_{n \to \infty}|a_n/b_n|=0$. 
      \item The notation $\mathbb{P}$ refers to a probability distribution on the probability space, where a random variable lives, and $\mathbb{E}$ denotes the expectation of the random variable.
    \item $\mu_A$ denotes a measure on a set $A$. If $A$ is a $n$-dimensional (sub)manifold, then $\Leb_A$ denotes the $n$-dimensional Lebesgue measure on $A$.
    \item $\mathbb{N}=\{1,2,3,\cdots\}$, $\mathbb{N}_0=\{0,1,2,3,\cdots\}$. 
\end{enumerate}

 %%%%%%%%%%%%%%%%%%%%%%%%%%%%%%%%%%%%%%%%%%%%%%%%%%%%%%%%%%%%%%%%%

\section{Ergodic theorems, Large deviations and Maximal large deviations}

We start with basic theorems in ergodic theory. 
\begin{definition}[Ergodic theorems]\ \par
    A fundamental Birkhoff's ergodic theorem says that if $\mu_{\mathcal{M}}$ is a probability measure, preserved by an ergodic dynamics $f:\mathcal{M} \to\mathcal{M}$, then almost surely for any function $\phi$ with $\int \phi d\mu_{\mathcal{M}}=0$,
\begin{equation*}\lim_{n\to \infty}\frac{\sum_{i \le n} \phi \circ f^i}{n}=0.  
\end{equation*}
\end{definition}

Equivalently, for any $\epsilon>0$, $\mu_\mathcal{M}$-a.s. $x \in \mathcal{M}$ there is $N=N_x>0$, such that 
\begin{equation*}
   \sup_{n\ge N} \Big|\frac{\sum_{i \le n}\phi \circ f^i(x)}{n}\Big|\le \epsilon.
\end{equation*}

Large deviations theory is concerned with probability of outliers in the convergence of averages of the values of observables along trajectories of dynamical systems, that is,
\begin{equation*}
    \mu_\mathcal{M}\Big(\Big| \frac{\sum_{i \le n}\phi \circ f^i}{n}\Big|\ge \epsilon\Big). 
\end{equation*}

 The paper \cite{vaientiadv} applied large deviations to obtain some geometric properties for certain smooth dynamical systems $f$, while in \cite{explicit}, the explicit estimates of constants in the upper bound of the error term in the large deviations principle were obtained.

If the dynamics is sufficiently strongly chaotic (e.g., if it has fast mixing rates or, in other words, correlations decay exponentially) and $\phi$ is sufficiently smooth, then there is a constant $C>0$, such that for any $n, N \ge 1$
\begin{equation*}
    \mu_\mathcal{M}\Big(\Big| \frac{\sum_{i \le n}\phi \circ f^i}{n}\Big|\ge \epsilon\Big) \approx e^{-Cn}, 
\end{equation*} and 
\begin{equation*}
    \mu_\mathcal{M}\Big(\sup_{n \ge N}\Big| \frac{\sum_{i \le n}\phi \circ f^i}{n}\Big|\ge \epsilon\Big) \approx \sum_{n\ge N}e^{-Cn} \approx e^{-CN}.
\end{equation*}

If the dynamics is weakly chaotic (i.e., slow mixing rates), then this estimate becomes weaker, \begin{equation*}
    \mu_\mathcal{M}\Big(\Big| \frac{\sum_{i \le n}\phi \circ f^i}{n}\Big|\ge \epsilon\Big) \precsim n^{-a} \text{ where }a>0.
\end{equation*}
(see \cite{melbourne09} for more details).

In this situation, in contrast with the case of the exponential estimate, a direct summation over $n$ does not give an analogous estimate 
\begin{equation*}
    \mu_\mathcal{M}\Big(\sup_{n \ge N}\Big| \frac{\sum_{i \le n}\phi \circ f^i}{n}\Big|\ge \epsilon\Big)\precsim N^{-a} \text{ where }a >0.
\end{equation*}

We call this type of estimate a \textbf{maximal large deviation}. 

One of the motivations to consider maximal large deviations is due to a rather typical approach (called an inducing method) in ergodic theory. This approach is usually aimed at finding a reference set $\Delta_0 \subsetneq \mathcal{M}$ and the first return time to it $R: \Delta_0 \to \mathbb{N}$, in order to prove some ``good" properties for the first return map $f^R: \Delta_0 \to \Delta_0$. Then this property is lifted to the original dynamics $f: \mathcal{M} \to \mathcal{M}$. In practice, this lifting technique relies on a quantitative estimate of the time $N$ in $\sup_{n\ge N} \frac{\sum_{i \le n}R \circ (f^R)^i}{n}$ (rather than $ \frac{\sum_{i \le N}R \circ (f^R)^i}{N}$). This is related to the presented below maximal large deviation approach.  

Let a transfer operators (Perron-Frobenius operator)  $P$ associated to $f$ is defined by the duality relations:
\begin{equation*}
    \int g \cdot P(h) d\mu_{\mathcal{M}}= \int g \circ f \cdot h d\mu_{\mathcal{M}} \text{ for all } h \in L^1(\mu_{\mathcal{M}})\text{, } g \in L^{\infty}(\mu_{\mathcal{M}}).
\end{equation*} 

Then we have the following.
\begin{theorem}[Maximal large deviations, see \cite{mldp}]\label{mld} \par
Suppose that $\phi\in L^{\infty}$ has zero mean, and there are constants $p\in \mathbb{N}, C_{\phi}>0, \beta \in (0, p)$ such that
\begin{equation*}
 ||P^n(\phi)||^p_p\le C_{\phi} n^{-\beta}  \text{ for any } n \ge 1.
\end{equation*} 

Then for any $\epsilon>0, N \ge 1$  \begin{equation*}
    \mu_\mathcal{M}\Big(\sup_{n \ge N}\Big| \frac{\sum_{i \le n}\phi \circ f^i}{n}\Big|\ge \epsilon\Big)\precsim_{\epsilon, p, \beta}\max\{||\phi||^{p}_{\infty}, C_{\phi}\}\cdot  ||\phi||^{p}_{\infty} \cdot N^{-\beta },
\end{equation*} and \begin{equation*}
    \Big|\Big|\sup_{n \ge N}\Big| \frac{\sum_{i \le n}\phi \circ f^i}{n}\Big|\Big|\Big|^{2p}_{2p}\precsim_{p, \beta}\max\{||\phi||^{p}_{\infty}, C_{\phi}\}\cdot  ||\phi||^{p}_{\infty} \cdot N^{-\beta}.
\end{equation*} 

In particular, if the function $\phi\in L^{\infty}$ satisfies the inequality $||P^n(\phi)||_1\le C_{\phi} n^{-\beta}$ for some $C_{\phi}>0$ and any $n \ge 1$, then for any $\epsilon>0$
     \begin{gather*}
     ||P^n(\phi)||^{1+\left \lceil \beta \right \rceil}_{1+\left \lceil \beta \right \rceil}\le C_{\phi} ||\phi||^{\left \lceil \beta \right \rceil}_{\infty} n^{-\beta},\\ 
     \mu_\mathcal{M}\Big(\sup_{n \ge N}\Big| \frac{\sum_{i \le n}\phi \circ f^i}{n}\Big|\ge \epsilon\Big)\precsim_{\epsilon,  \beta}\max\{||\phi||^{1+\left \lceil \beta \right \rceil}_{\infty}, C_{\phi} ||\phi||^{\left \lceil \beta \right \rceil}_{\infty}\}\cdot  ||\phi||^{1+\left \lceil \beta \right \rceil}_{\infty}  \cdot N^{-\beta }.
     \end{gather*}
\end{theorem}

Now we present applications to some concrete classes of dynamical systems. Denote by $d$ a Riemannian distance in $\mathcal{M}$, and by $\Leb$ a Riemannian volume form on $\mathcal{M}$. By $C^{1+}$ we denote the class of continuously differentiable maps with 
H\"older continuous derivatives. Let $f : \mathcal{M} \to \mathcal{M}$ be a piecewise $C^{1+}$-endomorphism of a Riemannian manifold $\mathcal{M}$.

\begin{definition}[Gibbs-Markov-Young (GMY) structures, see  \cite{Young2} or \cite{vaientiadv}]\label{GMY}\ \par 

Suppose that there exists a ball $\Delta_0\subset \mathcal{M}$,
a countable partition $\mathcal{P}$ (mod $0$) of $\Delta_0$ into topological balls $U$ with smooth boundaries, and a return time (not necessarily the first return time) function $R: \Delta_0 \to \mathbb{N}$, which is constant on each element of $\mathcal{P}$. We say that $f$ admits a Gibbs–Markov-Young (GMY) structure if the following properties hold:
\begin{enumerate}
    \item{Markov property:} for each $U\in \mathcal{P}$ and $R=R(U)$, $f^R: U\to \Delta_0$ is a $C^{1+}$-diffeomorphism (and, in particular, it is a bijection). Thus the induced map $F: \Delta_0 \to \Delta_0$, $F(x)=f^{R(x)}(x)$, is defined almost everywhere and satisfies a classical Markov property.
\item{Uniform expansion:} there exists $\rho\in (0,1)$ such that for almost all $x \in \Delta_0$ we have $||DF(x)^{-1}||\le \rho$. In particular the separation time $s(x,y)$, which is equal to the smallest integer $i\ge 0$ such that $F^i
(x)$ and $F^i
(y)$ lie in distinct elements of the partition $\mathcal{P}$, is defined and is  finite for almost all $x, y \in \Delta_0$.
\item{Bounded distortion:} there exists $C > 0$ such that for any two points $x,y \in \Delta_0$ with $s(x,y)\in [1,\infty) $ \begin{equation*}
    \Big|\frac{\det DF(x)}{\det DF(y)}-1\Big|\le C  \rho^{s(F(x),F(y))},
\end{equation*}where $\det DF$ is the Radon-Nikodym derivative of $F$ with respect to $\Leb$ on $\Delta_0$.
\end{enumerate}
\end{definition}

Theorem \ref{mld} implies the following corollary, see a proof in \cite{mldp}. 
\begin{corollary} \label{optimalmld}If $f$ admits  a Gibbs–Markov-Young structure, and $\Leb(x\in \Delta_0: R(x)>N)\approx N^{1+\beta}$ for some $\beta>0$, then there exists an invariant probability measure $\mu_{\mathcal{M}}$ and an open and dense set $H_0$ in the space of H\"older functions, such that for any zero mean function $\phi\in H_0$, and for any $\epsilon>0$,
\begin{equation*}
    \lim_{N\to \infty}\frac{\log \mu_{\mathcal{M}}\Big(\sup_{n \ge N}\Big| \frac{\sum_{i \le n}\phi \circ f^i}{n}\Big|\ge \epsilon\Big)}{\log N}=-\beta .
\end{equation*}
\end{corollary}

\section{Recurrences and Poisson Approximations}

The studies of Poisson approximations of the process of recurrences to small subsets in the phase spaces of chaotic dynamical systems are developed now into a large active area in the  dynamical systems theory. Another view at this type of problem is a subject of the theory of open dynamical systems where some positive measure subset $A$ of the phase space is named a hole, and hitting and escape through the hole processes are studied. The third view at this type of problems concerns statistics of extreme events (``record values") in the theory of random processes \cite{haydncmp, FFMA}. In this section we present some new advances in this area.

In a general set up, one picks a small measure subset $A$ in the phase space $\mathcal{M}$ of hyperbolic (chaotic) ergodic dynamical system and attempts to prove that in the limit, when the measure of $A$ approaches  zero, the corresponding process of recurrences to $A$ converges to the Poisson process. 

This area received an essential boost in L-S.Young papers \cite{Young1,Young2}, where a new general framework for analysis of statistical properties of hyperbolic dynamical systems was introduced. This approach employs representation of the phase space of a dynamical system as a tower (later called Young tower, Gibbs-Markov-Young tower, etc), which allow to study dynamics by analysing  recurrences to the base of this tower. Several developments of this method were proposed later, essentially all focused on the dynamical systems with weak hyperbolicity (slow decay of correlations). For such systems the method of inducing was employed, when the base of the tower is chosen as such subset of the phase space where the induced dynamics, generated by the recurrences to the base, is strongly hyperbolic.

The main results (Theorems \ref{thm2} and \ref{thm}) in this section are dealing with convergence of a random process, generated by the measure preserving dynamics, to the functional Poisson law in the total variation (TV) norm. We also obtain estimates of the corresponding convergence rates in the following form. For almost every $z\in \mathcal{M}$,\begin{align}\label{equaintro}
   d_{TV}\left(\mathcal{N}^{r,z,T}, Poisson \right) \precsim_{T,z} r^{a},
\end{align}where $Poisson$ is a Poisson point process and $\mathcal{N}^{r,z,T}$ is a dynamical point process which counts a number of entrances by an orbit to a metric ball $B_{r}(z)$ with radius $r$ and the center $z$ in the phase space of a dynamical system during a time interval $[0,T]$ (see their definitions below).  The notation $\precsim_{T,z}$ means that a constant in (\ref{equaintro}) depends only on $z$ and $T$. 

\begin{definition}[Dynamical point processes]\label{dynamicptprocess}\ \par
Let $(\mathcal{M},d)$ be a Riemannian manifold (with or without boundaries, connected or non-connected, compact or non-compact), $d$ is the Riemannian metric on $\mathcal{M}$ and $B_{r}(z)$ is a geodesic ball in $\mathcal{M}$ with a radius $r$ and a center $z \in \mathcal{M}$. We assume that dynamics $f: (\mathcal{M},\mu_{\mathcal{M}}) \to (\mathcal{M}, \mu_{\mathcal{M}})$ is ergodic with respect to (w.r.t.) some invariant probability measure $\mu_{\mathcal{M}}$.

Let $T>0$. Consider a dynamical point process on $[0,T]$, so that for any $t \in [0,T]$ 
\[\mathcal{N}^{r,T, z}_{t}:=\sum_{i=0}^{ t/\mu_{\mathcal{M}}(B_r(z))} \mathbbm{1}_{B_r(z)} \circ f^i.\]

Thus the dynamical point process $\mathcal{N}^{r,T, z}$ is a random counting measure on $[0,T]$.

\end{definition}

\begin{definition}[Poisson point processes]\label{poisson}\ \par
For any $T>0$, we say that $Poisson$ is a Poisson point process on $[0,T]$ if
\begin{enumerate}
    \item $Poisson$ is a random counting measure on $[0,T]$.
    \item $Poisson(A)$ is a Poisson-distributed random variable for any Borel set $A \subseteq [0,T]$.
    \item If $A_1, A_2, \cdots, A_n \subseteq [0,T]$ are pairwise disjoint, then  $Poisson(A_1), \cdots, Poisson(A_n)$ are independent.
    \item $\mathbb{E}Poisson(A)=\Leb(A)$ for any Borel set $A\subseteq [0,T]$.
\end{enumerate}

\end{definition}

\begin{definition}[Total variation norms of point processes]\label{totalnorm}\ \par
For any $T>0$ consider the $\sigma$-algebra $\mathcal{C}$ on the space of point processes on $[0,T]$, defined as \begin{align}\label{sigmaalg}
\sigma \left\{\pi^{-1}_AB: A \subseteq [0,T],B \subseteq \mathbb{N}\right\},
\end{align}
where $A, B$ are Borel sets and
$\pi_A$ is an evaluation map defined on the space of counting measures, so that for any counting measure $N$
\[\pi_AN:=N(A).\]

Now we can define the total variation norm for the Poisson approximation of a dynamical point process as
\[d_{TV}\left(\mathcal{N}^{r,T,z},Poisson\right):=\sup_{C\in \mathcal{C}}\left|\mu_{\mathcal{M}}(\mathcal{N}^{r,T,z} \in C)-\mathbb{P}(Poisson \in C)\right|\]
\end{definition}
\begin{definition}[Poisson approximations for small holes with convergence rates]\label{convergencerate}\ \par
Suppose that for any $T>0$ there exists a constant $a>0$ s.t. for almost every  $z\in \mathcal{M}$
\[d_{TV}\left(\mathcal{N}^{r,T,z},Poisson\right) \precsim_{T,z} r^a \to 0.\]

Then $a$ is called a convergence rate of a Poisson approximation.
\end{definition}
\begin{remark}
The convergence rates imply that for any sub-interval $[T_1, T_2]\subseteq [0,T]$, \[\sup_{C\subseteq \mathbb{N}}\left|\mu(\mathcal{N}^{r,T,z}[T_1,T_2]\in C)-\mathbb{P}(Poisson[T_1,T_2] \in C)\right|\precsim_{T,z} r^a \to 0.\]
\end{remark}

We now turn to hyperbolic Gibbs-Markov-Young structures \cite{Young1, Young2}:

\begin{definition}[Hyperbolic Gibbs-Markov-Young structures]\label{gibbs}\ \par

Introduce at first several notions concerning hyperbolic dynamics $f$ on Riemannian manifolds  $(\mathcal{M},d)$.

\begin{enumerate}
    \item An embedded disk $\gamma^u$ is called an unstable manifold if for every $x,y \in \gamma^u$ 
    \[\lim_{n \to \infty}d\Big(f^{-n}(x), f^{-n}(y)\Big)=0\]
    \item An embedded disk $\gamma^s$ is called a stable manifold if for every $x,y \in \gamma^s$ 
    \[\lim_{n \to \infty}d\Big(f^{n}(x), f^{n}(y)\Big)=0\]
    \item $\Gamma^u:=\{\gamma^u\}$ is called a continuous family of $C^1$-unstable manifolds if there is a compact
set $K^s$, a unit disk $D^u$ in some $\mathbb{R}^n$ and a map $\phi^u: K^s \times D^u \to \mathcal{M}$ such that 
\begin{enumerate}
    \item $\gamma^u=\phi^u\left(\{x\}\times D^u\right)$ is an unstable manifold,
    \item $\phi^u$ maps $K^s \times D^u$ homeomorphically onto its image,
    \item $x \to \phi^u|_{\{x\} \times D^u}$ defines a continuous map from $K^s$ to $\Emb^1\left(D^u, \mathcal{M}\right)$, where $\Emb^1\left(D^u, \mathcal{M}\right)$ is the space of $C^1$-embeddings of $D^u$ into $\mathcal{M}$.
\end{enumerate}
\end{enumerate}

A continuous family of $C^1$-stable manifolds $\Gamma^s:=\{\gamma^s\}$ is defined similarly.

\item We say that a compact set $ \Lambda \subseteq \mathcal{M}$ has a
hyperbolic product structure if there exist  continuous families of stable manifolds $\Gamma^s:=\{\gamma^s\}$ and of unstable manifolds $\Gamma^u:=\{\gamma^u\}$ such that

\begin{enumerate}

    \item $\Lambda=\left(\bigcup \gamma^s\right) \bigcap \left(\bigcup \gamma^u\right)$,
    
    \item $\dim \gamma^s+\dim \gamma^u=\dim \mathcal{M}$,
    
    \item each $\gamma^s$ intersects each $\gamma^u$ at
    exactly one point,
    
    \item stable and unstable manifolds are transversal, and the angles between them are uniformly bounded away from 0.
\end{enumerate}

A subset $\Lambda_1 \subseteq \Lambda$ is called a s-subset if $\Lambda_1$ has a hyperbolic product structure and, moreover, the corresponding families of stable and unstable manifolds  $\Gamma^s_1$ 
and $\Gamma^u_1$ can be chosen so that  $\Gamma^s_1 \subseteq \Gamma^s$ and $\Gamma^u_1 = \Gamma^u$. 

Analogously, a subset $\Lambda_2 \subseteq \Lambda$
is called an $u$-subset if $\Lambda_2$ has a hyperbolic product structure and the families $\Gamma^s_2$
and $\Gamma^u_2$ can be chosen so that  $\Gamma^u_2 \subseteq \Gamma^u$ and $\Gamma^s_2 = \Gamma^s$.

For $x \in \Lambda$, denote by $\gamma^u(x)$ (resp. $\gamma^s(x)$) the element of $\Gamma^u$ (resp. $\Gamma^s$)
which contains $x$. Also, for each $n \ge 1$, 
denote by $\left(f^n\right)^u$ the restriction of the map $f^n$ to $\gamma^u$-disks, and by $\det D\left(f^n\right)^u$ denote the Jacobian of $\left(f^n\right)^u$.

We say that the set $\Lambda$ with  hyperbolic product structure has also a \textbf{hyperbolic Gibbs-Markov-Young structure} if  the following  properties are satisfied 

\begin{enumerate}
    \item  Lebesgue detectability:  there exists $\gamma \in \Gamma^u$ such that $\Leb_{\gamma}\left(\Lambda \bigcap \gamma\right) > 0$.
    \item Markovian property: there exist pairwise disjoint $s$-subsets $\Lambda_1,\Lambda_2, \cdots \subseteq \Lambda$ such that
    \begin{enumerate}
        \item $\Leb_{\gamma}\left(\Lambda \setminus \left(\bigcup_{i \ge 1}\Lambda_i\right)\right)=0$ on each $\gamma \in \Gamma^u$,
        \item for each $i \ge 1$ there exists $R_i \in \mathbb{N}$ such that $f^{R_i} (\Lambda_i)$ is an $u$-subset, and for all $x \in \Lambda_i$
        \[f^{R_i}\left(\gamma^s(x)\right) \subseteq \gamma^s\left(f^{R_i}(x)\right)\]
        and
    \[f^{R_i}\left(\gamma^u(x)\right) \supseteq \gamma^u(f^{R_i}(x)).\]
    \end{enumerate}
   
Define now a return time function $R : \Lambda \to \mathbb{N}$ and a return function $f^R: \Lambda \to \Lambda$, so that for each $i \ge 1$ 
\[R\big|_{\Lambda_i}=R_i
\text{ and }
f^R\big|_{\Lambda_i}=f^{R_i}\big|_{\Lambda_i}\]

Next, the separation time $s(x, y)$ for $x, y \in \Lambda$ is defined as
\[s(x,y):=\min \{n \ge 0: \left(f^R\right)^n(x) \text{ and } \left(f^R\right)^n(y) \text{ belong to different sets } \Lambda_i\}.\]

We also assume that there are constants $C > 1, \alpha > 0$ and $0 < \beta < 1$, which depend only on $f$ and $\Lambda$, such that the following conditions hold

\item Polynomial contraction on stable leaves. For any $\gamma^s \in \Gamma^s, x,y \in \gamma^s,  n \ge 1$,
\[ d\Big(f^n(x), f^n(y)\Big) \le Cn^{-\alpha}.\]

\item Backward polynomial contraction on unstable leaves. For any $ \gamma^u \in \Gamma^u, x,y \in \gamma^u,  n\ge 1$,
\[ d\Big(f^{-n}(x), f^{-n}(y)\Big) \le Cn^{-\alpha}.\]
\item Bounded distortion: for any $ \gamma \in \Gamma^u \text{ and } x, y \in \gamma \bigcap \Lambda_i$ for some $\Lambda_i$, 
\[\log \frac{\det D\left(f^R\right)^u(x)}{\det D\left(f^R\right)^u(y)} \le C  \beta^{s\left(f^R(x), f^R(y)\right)}.\]
\item Regularity of the stable foliation. For each $\gamma, \gamma'\in \Gamma^u$ denote 
\[\Theta_{\gamma, \gamma'}: \gamma' \bigcap \Lambda \to \gamma \bigcap \Lambda: x \to \gamma^s(x)\bigcap \gamma.\]

Then the following properties hold
\begin{enumerate}
    \item $\Theta_{\gamma, \gamma'}$ is absolutely continuous and for any $ x \in \gamma \bigcap \Lambda$
    \[\frac{d\left(\Theta_{\gamma, \gamma'}\right)_{*}\Leb_{\gamma'}}{d\Leb_{\gamma}}(x)=\prod_{n \ge 0} \frac{\det Df^u\left(f^n\left(x\right)\right)}{\det Df^u\left(f^n\left(\Theta_{\gamma, \gamma'}^{-1}x\right)\right)},\]
    \[ \frac{d\left(\Theta_{\gamma, \gamma'}\right)_{*}\Leb_{\gamma'}}{d\Leb_{\gamma}}(x)= C^{\pm 1},\]
    \item for any $x,y \in \gamma \bigcap \Lambda$ 
    \[\log \frac{\frac{d\left(\Theta_{\gamma, \gamma'}\right)_{*}\Leb_{\gamma'}}{d\Leb_{\gamma}}(x)}{\frac{d\left(\Theta_{\gamma, \gamma'}\right)_{*}\Leb_{\gamma'}}{d\Leb_{\gamma}}(y)} \le C  \beta^{s(x,y)}.\]
    
\end{enumerate}
\item Aperiodicity: $\gcd\left(R_i, i \ge 1\right)= 1$. 
   \item Decay rate of the return times $R$: there exist $\xi>1$ and $\gamma\in \Gamma^u$ such that \begin{align*}
       \Leb_{\gamma}\left(R>n\right) \le Cn^{-\xi}.
   \end{align*}
\end{enumerate}

\textbf{SRB measures}: 
Let the dynamics $f: \mathcal{M} \to \mathcal{M}$ has hyperbolic Gibbs-Markov-Young structure. It was proved in \cite{Young2, Young1} that there exists  an ergodic probability measure $\mu_{\mathcal{M}}$  such that for any unstable manifold $\gamma^u$ (including $\Gamma^u$) we have $\mu_{\gamma^u} \ll \Leb_{\gamma^u}$, where $\mu_{\gamma^u}$ is the conditional measure of $\mu_{\mathcal{M}}$ on an unstable manifold $\gamma^u$. Such measure $\mu_{\mathcal{M}}$ is called Sinai-Ruelle-Bowen measure (SRB measure).

\end{definition}

\begin{assumption}[Geometric regularities]\label{geoassumption}\ \par
Assume that $f: \mathcal{M}\to \mathcal{M}$ has the hyperbolic Gibbs-Markov-Young structure, as described in Definition \ref{gibbs}, and
\begin{enumerate}
    
   \item $f$ is bijective and a local $C^1$-diffeomorphism on $\bigcup_{i \ge 1} \bigcup_{0 \le j < R_i} f^j(\Lambda_i)$.
   \item The following limit 

  \[ \dim_H\mu_{\mathcal{M}}:=\lim_{r \to 0} \frac{\log \mu_{\mathcal{M}}(B_r(z))}{\log r}\]
   
exists  for almost every  $z\in \mathcal{M}$. Then  $dim_H\mu_{\mathcal{M}}$ is  called a Hausdorff dimension of the measure $\mu_{\mathcal{M}}$.
 
 \item $\alpha \dim_H\mu_{\mathcal{M}}>1$, where $\alpha$ is defined in Definition \ref{gibbs}.
\end{enumerate}
\end{assumption}

\begin{assumption}[The first returns \& interior assumptions on $\Lambda$]\label{assumption}\ \par

Assume that $f: \mathcal{M}\to \mathcal{M}$ has the hyperbolic Gibbs-Markov-Young structure, and there are constants $C>1$ and $\beta\in (0,1)$ (the same as in Definition \ref{gibbs}) such that
    
\begin{enumerate}   
    \item $R: \Lambda \to \mathbb{N}$ is the first return time and $f^R: \Lambda \to \Lambda$ is the first return map for $\Lambda$.
    
    This implies that $f^R$ is bijective.
    \item  for any $\gamma \in \Gamma^s, \gamma_1 \in \Gamma^u, x, y \in \gamma \bigcap \Lambda, x_1, y_1 \in \gamma_1 \bigcap \Lambda $,
    \[d\Big(\left(f^R\right)^n(x), \left(f^R\right)^n(y)\Big) \le C \beta^n,\]
    and
    \[d\Big(\left(f^R\right)^{-n}(x_1), \left(f^R\right)^{-n}(y_1)\Big) \le C \beta^n  d(x_1,y_1).\]

    \item $\mu_\mathcal{M}\{\interior{(\Lambda)}\}>0$ and $\mu_\mathcal{M}(\partial \Lambda)=0$, where \begin{align*}
        \interior{\Lambda}:=\left\{x \in \Lambda: \text{ there exists } r_x>0 \text{ s.t. } \mu_\mathcal{M}\left(B_{r_x}(x)\setminus \Lambda\right)=0 \right\},\quad \partial \Lambda:=\Lambda \setminus \interior{\Lambda}. 
    \end{align*}    
    
    In other words, $x \in \interior{\Lambda}$ if and only if $x \in \Lambda$ and there is a small ball $B_{r_x}(x)$ s.t. $B_{r_x}(x) \subseteq \Lambda$ $\mu_\mathcal{M}$-almost surely. 
\end{enumerate}

\end{assumption}

We formulate now the first result in this section.

\begin{theorem}[Convergence rates for Poisson laws \Romannum{1}, see \cite{Sucmp}]\label{thm1}\ \par
Assume that the dynamics $f: (\mathcal{M}, \mu_{\mathcal{M}}) \to (\mathcal{M}, \mu_{\mathcal{M}})$ has a first return hyperbolic Gibbs-Markov-Young structure (see Definition \ref{gibbs}) and satisfies Assumptions \ref{geoassumption} and \ref{assumption}. Then for any $T>0$ the following results hold 

\begin{enumerate}
    \item $\dim_H\mu_\mathcal{M} \ge \dim {\gamma}^u$ 
    \item if either $\alpha> \frac{2}{\dim \gamma^u}-\frac{1}{\dim_H \mu}$ or  $\mu_{\mathcal{M}} \ll \Leb_{\mathcal{M}}$ with $\frac{d\mu_\mathcal{M}}{d\Leb_{\mathcal{M}}} \in L^{\infty}_{loc}(\mathcal{M})$, then for almost every $z \in \mathcal{M}$
    \[d_{TV}\left(\mathcal{N}^{r,z,T}, Poisson\right) \precsim_{T,\xi,z} r^{a},\]
    where a computable constant $a>0$ depends on $\xi>1$, $\dim_H \mu_{\mathcal{M}}, \dim \gamma^u $ and $\alpha$, but it does not depend on $z\in \mathcal{M}$. 
    
\end{enumerate}

\end{theorem}

\begin{definition}[Induced measurable partitions]\label{mpartition}\ \par
We say that a probability measure $\mu_\mathcal{M}$ for the dynamics $f: \mathcal{M} \to \mathcal{M}$ has an induced measurable partition if there are constants $\beta \in (0,1), C>1$ (the same as in the Definition \ref{gibbs}) and $b>0$ such that
\begin{enumerate}
    \item There exists a subset $X \subseteq \mathcal{M}$ with $\mu_{\mathcal{M}} \{\interior{(X)}\}>0$, $\mu_{\mathcal{M}}(\partial X)=0$.
    \item The subset $X$ has a measurable  partition $\Theta:=\{\gamma^u(x)\}_{x \in X}$ (which could be different from $\Gamma^u$), such that the elements of $\Theta$ are  disjoint  connected unstable manifolds, so that $\mu_{\mathcal{M}}$-almost surely $X=\bigsqcup_{x \in X} \gamma^u(x)$ and for any function $g$
    \[\mu_{X}(g)=\int_X \mu_{\gamma^u(x)}(g) d\mu_X(x),\]
    where $\mu_X:=\frac{\mu_\mathcal{M}|_{X}}{\mu_\mathcal{M}(X)}$ and $\mu_{\gamma^u(x)}$ is the conditional probability induced by $\mu_{X}$ on $\gamma^u(x)\in \Theta$.
    \item Each $\gamma^u\in \Theta$ is (at least $C^1$) smooth.
    \item All $\gamma^u\in \Theta$ have uniformly bounded sectional curvatures and the same dimensions.
    \item For any $\epsilon \in (0,1)$\[\mu_X \{x \in X: |\gamma^u(x)|< \epsilon\}\le C  \epsilon^b,\] where $|\gamma^u(x)|$ is the radius of the largest inscribed geodesic ball  in $\gamma^u(x) \in \Theta$, where a 
    geodesic ball is defined with respect to the distance   $d_{\gamma^u(x)}$ on $\gamma^u(x)$, induced by the Riemannian metric. This property implies that almost every $\gamma^u(x)\in \Theta$ is non-degenerate, i.e., $|\gamma^u(x)|>0$ for almost every $x\in X$.
    \item For almost every point $x\in X$ we have $\mu_{\gamma^u(x)} \ll \Leb_{\gamma^u(x)}$, $\mu_{\gamma^u(x)}\left(\gamma^u(x)\right)>0$, and for any $y, z \in\gamma^u(x)$
   \[\frac{d\mu_{\gamma^u(x)}}{d\Leb_{\gamma^u(x)}}(y)=C^{\pm 1}  \frac{d\mu_{\gamma^u(x)}}{d\Leb_{\gamma^u(x)}}(z).\]
    \item Denote by $\overline{R}$ the first return time to $X$ for $f$. Then the first return map $f^{\overline{R}}: X \to X$ has an exponential u-contraction, i.e., for any $\gamma^u \in \Theta, x,y \in \gamma^u, n \ge 1$
    \[d\Big(\left(f^{\overline{R}}\right)^{-n}(x),\left(f^{\overline{R}}\right)^{-n}(y)\Big) \le C \beta ^n d(x,y),\]
    and an exponential decay of correlations, i.e., for any $h \in \Lip(X)$\begin{align*}
        \left|\int h \circ \left(f^{\overline{R}}\right)^n  h d\mu_{X}-\left(\int h d\mu_X\right)^2\right|\le C\beta^n ||h||^2_{\Lip}. 
    \end{align*}
\end{enumerate}
\end{definition}

Now we can formulate the second result in this section. 

\begin{theorem}[Convergence rates for the Poisson laws  \Romannum{2}, \cite{Sucmp}]\label{thm2}\ \par
Assume that the dynamics $f: (\mathcal{M}, \mu_{\mathcal{M}}) \to (\mathcal{M}, \mu_{\mathcal{M}})$ has the hyperbolic Gibbs-Markov-Young structure (see Definition \ref{gibbs}), satisfies the Assumption \ref{geoassumption}, and $\mu_\mathcal{M}$ has an induced measurable partition (see Definition \ref{mpartition}). Then for any $T>0$, the following results hold.

\begin{enumerate}
    \item $\dim_H \mu_{\mathcal{M}} \ge \frac{ b}{b+\dim \gamma^u} \dim \gamma^u$ and 
    \item if either $\alpha>  \frac{2}{\dim \gamma^u} \frac{b+\dim \gamma^u}{b}-\frac{1}{\dim_H \mu_{\mathcal{M}}} $ or  $\mu_{\mathcal{M}} \ll \Leb_{\mathcal{M}}$ and $\frac{d\mu_{\mathcal{M}}}{d\Leb_{\mathcal{M}}} \in L^{\infty}_{loc}(\mathcal{M})$, then for almost every (a.e.) $z \in \mathcal{M}$,
    \[d_{TV}\left(\mathcal{N}^{r,z,T}, Poisson\right) \precsim_{T, \xi, z,b} r^{a},\]
    where a constant $a>0$ depends on $\xi>1$, $\dim_H \mu_{\mathcal{M}}, \dim \gamma^u, b$ and $\alpha$, but it does not depend on $z\in \mathcal{M}$. 
\end{enumerate}

\end{theorem}

We will consider now some applications of the above Theorems.
\subsection{Intermittent solenoids}\label{appsolenoid}
Following \cite{Alves} let $\mathcal{M}=S^1 \times \mathbb{D}$, $f_{\gamma}(x,z)=(g_{\gamma}(x), \theta  z+e^{2\pi i x}/2)$, where $g_{\gamma}: S^1 \to S^1$ is a continuous map of degree $d\ge 2$ and $\gamma \in (0, +\infty)$ such that 
\begin{enumerate}
    \item $g_{\gamma}$ is $C^2$ on $S^1 \setminus \{0\}$ and $Dg_{\gamma}>1$ on $S^1 \setminus \{0\}$,
    \item $g_{\gamma}(0)=0, Dg_{\gamma}(0+)=1$ and  $xD^2g_{\gamma}(x) \sim x^{\gamma}$ for sufficiently small positive $x$,
    \item $Dg_{\gamma}(0-)>1$,
    \item $\theta>0$ is so small that $\theta  ||Dg_{\gamma}||_{\infty}< 1-\theta$.
\end{enumerate}

By Theorem \ref{thm1} we have the following
\begin{corollary}\label{fplsolenoid}
Poisson limit laws  hold for $f_{\gamma}$ for any $\gamma \in (0,1)$ with convergence rates  (\ref{equaintro}) . 
\end{corollary}

\subsection{Axiom A attractors}\label{appA}

It was well-known that for Axiom A attractors $\Sigma \subsetneq \mathcal{M}$ with $\dim \gamma^u=1$ the Poisson limit laws hold. We will show that conditions on $\dim_H\mu$ and $\dim \gamma^u$ can be dropped. 

\begin{definition}[Axiom A attractors]\ \par

Let $f : \mathcal{M} \to  \mathcal{M}$ be a  $C^2$-diffeomorphism. A compact set $\Sigma \subseteq \mathcal{M}$ is called an Axiom A attractor
if 
\begin{enumerate}
    \item There is a neighborhood $U$ of $\Sigma$, called its basin, such that $f^n(x) \to \Sigma$ for every $x \in U$.
\item The tangent bundle over $\Sigma$ is split into $E^u \oplus E^s$,  where $E^u$ and $E^s$ are $df$-invariant subspaces. 
\item $df|_{E^u}$ is uniformly expanding, and $df|_{E^s}$ is uniformly contracting.
\item $f: \Sigma \to \Sigma$ is topologically mixing.
\end{enumerate}
\end{definition}

The support of the SRB measure $\mu_{\mathcal{M}}$ is $\Sigma$. Applying now Theorem \ref{thm1} we obtain the following.

\begin{corollary}\label{axiomfpl}
Poisson limit laws hold for Axiom A attractors with convergence rates  (\ref{equaintro}).
\end{corollary}

\subsection{H\'enon attractors}\label{henon}
 
The Poisson limit laws for certain H\'enon attractors (see \cite{Young1}) have been proved in \cite{collet}. However, convergence rate for this class of dynamical systems was obtained  in a weaker form. A stronger rate of convergence  was derived (see \ref{equaintro}) by using Theorem \ref{thm2}.

\begin{corollary}\label{fplhenon}
Poisson limit laws with convergence rates (\ref{equaintro}) hold for H\'enon attractors which can be modelled by Young towers. 
\end{corollary}

\subsection{Billiards}
\begin{definition}[Billiard tables, billiard maps and phase spaces]\label{billiarddef} \par
We consider a billiard in a two-dimensional region $Q\subseteq \mathbb{R}^2$ (called a billiard table) with a piece-wise smooth (of class $C^3$) boundary $\partial Q$. Each smooth piece has a uniformly bounded curvature. The boundary $\partial Q$ is equipped with a field of inward unit normal vectors $n(q), q \in \partial Q$. 

A billiard is a dynamical system generated by the motion of a point
particle with the unit velocity inside the region $Q$, and being reflected from its
boundary according to the law ``the angle of incidence equals the angle of
reflection". It means that upon reflection the tangent component of the
velocity remains the same, while the normal component changes its sign
according to the rule $v_{+}=v_{-}-2\langle n(q), v_{-}\rangle n(q)$, where $v_{+}$ (resp. $v_{-}$) is the velocity of the particle immediately after (resp. before) reflection.

The phase space of a billiard is the restriction of the unit tangent
bundle of $\mathbb{R}^2$ to $Q$. We will use the standard notation for phase points $x=(q, v)$, where $q$ is the point of the configuration space $Q$ and $v$ is the unit velocity vector. The billiard preserves the Liouville measure $d\nu:=dq dv$ where $dq$ and $dv$ are Lebesgue measures on $Q$ and on the unit
one-dimensional sphere respectively. The corresponding flow will be denoted by $\{S^t\}$. It is customary for billiard-type systems to study instead of $\{S^t\}$ a
dynamical system with discrete time, which is called a billiard map $f$. Denote \[\mathcal{M}=\{x=(q, v), q \in \partial Q, \langle v, n(q)\rangle >0\}.\]  For $x=(q, v) \in \mathcal{M}$ let $\tau (x)$ be
the first positive moment of reflection of the billiard
orbit defined by $x$ off the boundary. Then the billiard map is defined by $f(x)=(q', v')=S^{\tau}
x$, so that $q'$ is the point of
the next reflection and $v'$ is the outcoming velocity vector at that point. We call $\mathcal{M}$ a phase space of the billiard map $f$.

Due to the structure of $\partial Q$,  the set of singular points of the boundary $\partial Q$ has measure zero. The angle $\phi$ of the velocity vector $v$ varies from $-\pi/2$ to $\pi/2$ at any regular point $q \in \partial Q$.  Hence $\mathcal{M}=\partial Q \times [-\pi/2,\pi/2]$ almost surely. In what follows we always identify $\mathcal{M}$ with $\partial Q \times [-\pi/2,\pi/2]$, and denote the phase point $x\in \mathcal{M}$ by $(q, \phi) \in \partial Q \times [-\pi/2,\pi/2]$ throughout this subsection.

The phase space  $\mathcal{M}=\partial Q \times [-\pi/2,\pi/2]$ is endowed with a natural Riemannian metric $d_{\mathcal{M}}$ and Riemannian volume $\Leb_{\mathcal{M}}$. The billiard map $f$ preserves the measure
\[d\mu_{\mathcal{M}}:=(2\Leb_{\partial Q} \partial Q)^{-1}\cos \phi d\phi dq=(2\Leb_{\partial Q} \partial Q)^{-1}\cos \phi \Leb_{\mathcal{M}},\]where $dq$ is the one-dimensional Lebesgue measure on the boundary $\partial Q$ and
$d\phi$ is the one-dimensional Lebesgue (uniform) measure on $[-\pi/2,\pi/2]$. 
\end{definition}

Before presenting the results, we introduce some more definitions.

\begin{definition}[An induced system]\label{inducesystem}\ \par
Suppose that there is a fixed subset $X\subseteq \mathcal{M}$ with $\Leb_{\mathcal{M}}(X)>0$. The first return time to $X$ is denoted as $R: X\to \mathbb{N}$. We assume that $X$ can be partitioned into countably many connected pieces
\begin{equation}\label{parition}
    X=\bigcup_{i\ge 1} X_i,
\end{equation} so that $R$ is constant on each $X_i$, and  \[\Leb_{\mathcal{M}} (\partial X_i)=0,\quad \interior{X_i} \bigcap \interior{X_j}=\emptyset \text{ for } i \neq j.\]

The first return time $R$ induces the first return map  $f^R: X \to X$ and a new dynamical system $(X, f^R)$.

\end{definition}
\begin{definition}[Singularities and (un)stable manifolds]\label{unstablemfd} \par
 Denote by $\mathbb{S}\subseteq X$ the singularity set for $f^R$. The set $\mathbb{S}$ has zero Lebesgue measure, and $\mathbb{S}^c \subseteq X$ consists of countably many open connected components. Unstable (resp. stable) manifolds are connected components of $(\bigcup_{i \ge 0}(f^R)^{i}\mathbb{S})^c$ (resp. $(\bigcup_{i \ge 0}(f^R)^{-i}\mathbb{S})^c$). A closed and connected part of the unstable (resp. stable) manifold will be called an unstable (resp. stable ) disk. We denote each unstable (resp. stable) manifold/disk by $\gamma^u$ (resp. $\gamma^s$), and its tangent vectors by $v^u$ (resp. $v^s$).
\end{definition}

\begin{remark}\label{singular}
The singularity set $\mathbb{S}$ consists of the points in $X $ which are not ``well-behaved". It includes the discontinuities and the points where the map $f^R$ is not differentiable. (It may also include other points in $X$ with some ``bad" properties).
\end{remark}
\begin{definition}[Chernov-Markarian-Zhang (CMZ) structures]\label{cmz}\ \par
We say that an induced system $(X, f^R)$ in Definition \ref{inducesystem} is a CMZ structure on $(\mathcal{M}, f)$ if there are constants $C>0$ and $\beta \in (0,1)$ such that the following conditions hold
\begin{enumerate}
    \item Hyperbolicity. For any $n\in \mathbb{N}$, $v^u$ and $v^s$,\begin{align*}
        |D(f^R)^n v^u|\ge C\beta^{-n} |v^u|, \quad  |D(f^R)^{n} v^s|\le C\beta^n |v^s|,
    \end{align*}where $|\cdot|$ is the Riemannian metric induced from $\mathcal{M}$ to (un)stable manifolds.
    \item SRB measures and u-SRB measures. The conditional probability measure $\mu_X:=\frac{\mu_{\mathcal{M}}|_X}{\mu_{\mathcal{M}}(X)} \ll \Leb_{X}$ with $\frac{d\mu_X}{d\Leb_X}\in L^{\infty}$ and $\frac{d\mu_X}{d\Leb_X}>0$ $\Leb_X$-a.s. on $X$. Moreover, conditional distributions on $\gamma^u$ (say $\mu_{\gamma^u}$) are absolutely continuous w.r.t. Lebesgue measure $\Leb_{\gamma^u}$ on $\gamma^u$.
    \item Distortion bounds. Let $d_{\gamma^u}(\cdot, \cdot)$ be the distance measured along $\gamma^u$. By $\det D^uf^R$ we denote the Jacobian of $Df^R$ along the unstable manifolds. Then, if $x, y\in X$ belong to a $\gamma^u$, and $(f^R)^n$ is $C^2$-smooth on $\gamma^u$, then the following relation holds
    \begin{align*}
        \log \frac{\det D^u(f^R)^n(x)}{\det D^u(f^R)^n(y)}\le \psi\left[d_{\gamma^u}\Big((f^R)^nx,(f^R)^ny\Big)\right],
    \end{align*}where $\psi(\cdot)$ is some function, which does not depend on $\gamma^u$, and $\lim_{s\to 0^+}\psi(s)=0$.
    \item Bounded curvatures. The curvatures of all $\gamma^u$ are uniformly
bounded by $C$.
\item Absolute continuity. Consider a holonomy map $h: \gamma^u_1 \to \gamma^u_2$,  which maps a point $x \in \gamma^u_1$ to the point $h(x)\in \gamma^u_2$, such that both $x$ and $h(x)$ belong to the same $\gamma^s$. We assume that the holonomy map satisfies the following relation \begin{align*}
    \frac{\det D^u(f^R)^n(x)}{\det D^u(f^R)^n\big(h(x)\big)}=C^{\pm 1} \text{ for all }n \ge 1 \text{ and }x\in \gamma^u_1,
\end{align*}
\item Growth lemmas. For any
sufficiently small $\delta>0$ and for any disk on a $C^1$-smooth unstable manifold  $\gamma^u$ with $\diam \gamma^u \le \delta_0$ denote by $U_{\delta} \subseteq \gamma^u$ a $\delta$-neighborhood of the subset $\gamma^u \bigcap \bigcup_{0\le i \le N}(f^R)^{-i}\mathbb{S}$ within the set $\gamma^u$. Then there exist $N \in \mathbb{N}$, sufficiently small $\delta_0>0$ and
constants  $ \kappa, \sigma > 0$ which satisfy the following condition. There exists an open subset $V_{\delta} \subseteq \gamma^u \setminus U_{\delta}$, such that $\Leb_{\gamma^u}\big(\gamma^u \setminus (U_{\delta}\bigcup V_{\delta})\big)=0$, and for any $\epsilon >0$ \begin{gather*}
    \Leb_{\gamma^u}(r_{V_{\delta},N}<\epsilon)\le 2\epsilon \beta+\epsilon C \delta_0^{-1} \Leb_{\gamma^u} (\gamma^u),\\
    \Leb_{\gamma^u}(r_{U_{\delta},0}< \epsilon)\le C \delta^{-\kappa} \epsilon,\\
    \Leb_{\gamma^u}(U_{\delta})\le C \delta^{\sigma},
\end{gather*}where $r_{U_{\delta},0}(x):=d_{\gamma^u}(x, \partial U_{\delta})$, $r_{V_{\delta}, N}(x):=d_{(f^R)^N\gamma^u}\big((f^R)^Nx, \partial (f^R)^NV_{\delta}(x)\big)$, and $V_{\delta}(x)$ is a connected component of $V_{\delta}$, which contains $x$.

\item Finiteness: $\int R d\mu_X< \infty$.
\item Mixing: $\gcd \{R\}=1$.

\end{enumerate}

\end{definition}

\begin{remark}
 In the growth lemmas a positive integer $N$ is usually chosen as a sufficiently large number.The conditions that $\gcd{R}=1$ and K-mixing of $f^R$ guarantee that $f$ is also K-mixing.
\end{remark}

Consider now the first return tower \begin{align*}
    \Delta:=\{(x,n)\in X \times \{0,1,2,\cdots\}: n < R(x)\}.
\end{align*}

Dynamics $F: \Delta \to \Delta$ is defined as $F(x, n)=(x, n+1)$ if $n+1\le R(x)-1$ and $F(x,n)=(f^Rx, 0)$ if $n=R(x)-1$. The projection $\pi:\Delta \to \mathcal{M}$ is defined by $\pi(x,n):=f^n(x)$ as $\pi \circ F=f \circ \pi$. Finally we introduce projections $\pi_{X}: \Delta \to X$ and $\pi_{\mathbb{N}}: \Delta \to \mathbb{N}_0$, so that for any $(x,n) \in \Delta$ \begin{equation}\label{projections}
    \pi_{X}(x,n)=x,\quad \pi_{\mathbb{N}}(x,n)=n.
\end{equation}

Extend now $\mu_X$ from $X$ to $\Delta$ as \begin{align*}
    \mu_{\Delta}:=(\int R d\mu_X)^{-1}\sum_{j}(F^j)_{*}\big(\mu_X|_{\{R>j\}}\big),
\end{align*} which reproduces the invariant probability measure on $\mathcal{M}$ \begin{align*}
\mu_{\mathcal{M}}=(\pi)_{*}\mu_{\Delta}.  \end{align*}

We identify $\Delta_0:=(X \times \{0\})\bigcap \Delta$ with $X$, $\mu_{\Delta_0}$ with $\mu_X$ and $F^R$ with $f^R$. Therefore $\pi: X \to X$ is the identity map.

Note that $\pi:\Delta \to \mathcal{M}$ is bijective. Thus $(\Delta,F)$ is identical to $(\mathcal{M},f)$, and $(X, f^R)=(\Delta_0, F^R)$ is a CMZ structure on $(\Delta, F)$.

\begin{remark}\ \par
\begin{enumerate}
    \item If $R=1$, then $X=\mathcal{M}$, and thus $(\mathcal{M}, f, \mu_{\mathcal{M}})$ has a CMZ structure.
    \item It follows from \cite{Chernovjsp, CZnon} that $(X, f^R, \mu_X)$ can be modelled by a hyperbolic Young tower \cite{Young1}.
    \item $(\mathcal{M}, f, \mu_{\mathcal{M}})$ is a K-system (and therefore mixing) because of the condition that $f^R$ is K-mixing and $\gcd \{R\}=1$.
\end{enumerate}
\end{remark}
\begin{definition}[Holes and dynamical point processes]\label{dynamicalpointprocess} \par
Throughout this subsection, a hole within the boundary $\partial Q$ is an open disk $B_{r}(q)$ with radius $r$ and the center at a regular point $q \in \partial Q$ of the boundary of a billiard table. 

Similar to Definition \ref{dynamicptprocess}, define now a different dynamical point process $\mathcal{N}^{r,q}$ on $\mathbb{R}^{+}\bigcup\{0\}$. For any measurable $A\subseteq \mathbb{R}^{+}\bigcup\{0\}$, and any $x \in \mathcal{M}$. \begin{align*}
    \mathcal{N}^{r,q}(x)(A):&=\#\{i\ge 0: f^i(x) \in B_r(q) \times [-\pi/2,\pi/2], \quad i \cdot \mu_{\mathcal{M}}\big(B_r(q) \times [-\pi/2,\pi/2]\big) \in A\}\\
    &=\sum_{i \cdot \mu_{\mathcal{M}}(B_r(q)\times [-\pi/2,\pi/2])\in A}\mathbbm{1}_{B_r(q)\times [-\pi/2,\pi/2]}\circ f^i(x).
\end{align*}

We will usually drop the symbol $x$ and write 
 \begin{align*}
    \mathcal{N}^{r,q}(A)=\sum_{i \cdot \mu_{\mathcal{M}}(B_r(q)\times [-\pi/2,\pi/2])\in A}\mathbbm{1}_{B_r(q)\times [-\pi/2,\pi/2]}\circ f^i.
\end{align*}

By using $\mu_{\mathcal{M}}:=(\pi)_{*}\mu_{\Delta}$, we get 
 \begin{align*}
    \mathcal{N}^{r,q}(A)=\sum_{i \cdot \mu_{\Delta}(A_r)\in A}\mathbbm{1}_{A_r}\circ F^i,
\end{align*}where $A_{r}:=\pi^{-1}(B_r(q)\times [-\pi/2,\pi/2])$.
\end{definition}
\begin{remark}
Following the theory of point processes  we have that 
\[\mathcal{N}^{r,q}(x)=\sum_{i\ge 0: f^i(x) \in B_r(q) \times [-\pi/2,\pi/2]} \delta_{i\cdot \mu_{\mathcal{M}}(B_r(q)\times [-\pi/2, \pi/2]) },\] where $\delta$ is a Dirac measure. Hence, $\mathcal{N}^{r,q}$ is a random counting measure, e.g., $\mathcal{N}^{r,q}(x)[0,1]$ counts the number of $i \in [0, \frac{1}{\mu_{\mathcal{M}}\{B_r(q)\times [-\pi/2, \pi/2]\}}]$, such that $f^i(x)$ belongs to $B_r(q)\times [-\pi/2, \pi/2]$. $B_r(q)\times [-\pi/2, \pi/2]$ is a large hole in the phase space $\mathcal{M}$.
\end{remark}

\begin{definition}[Poisson approximations for large holes]\label{pdef} \par
We say that $\mathcal{N}^{r,q} \to_d Poisson$ (see Definition \ref{poisson}) if  $\mathcal{N}^{r,q}g \to_d Poisson(g)$ for any $g\in C^{+}_c(\mathbb{R}^+\bigcup\{0\})$, i.e.,
\begin{align*}
    \lim_{r\to 0}\int \exp{(-t\cdot  \mathcal{N}^{r,q}g)}d\mu_{\mathcal{M}}=\int \exp{[-t\cdot Poisson(g)]}d\mathbb{P} \text{ for all }t>0,
\end{align*} where $C^{+}_c(\mathbb{R}^+\bigcup \{0\})$ is the space consisting of positive continuous functions with a compact support, defined on $[0,\infty)$.

This is equivalent to the relation \begin{align*}
    (\mathcal{N}^{r,q}I_1,\mathcal{N}^{r,q}I_2, \cdots, \mathcal{N}^{r,q}I_k) \to_d (Poisson(I_1), Poisson(I_2), \cdots, Poisson(I_k)) 
\end{align*}for any $k \in \mathbb{N}$ and any bounded intervals $I_1, \cdots I_k\subseteq{R}^+\bigcup\{0\}$.

In other words, the limit distribution of $\mathcal{N}^{r,q}$ is Poissonian, when the disk $B_r(q)$ shrinks to a regular point on the boundary $\partial Q$ of a billiard table.

\end{definition}

\begin{remark}
    If $I_1,\cdots, I_k$ are restricted inside an interval $[0,T]$, the convergence in this definition is equivalent to the  $d_{TV}$-convergence in Definition \ref{convergencerate} without convergence rates.
\end{remark}

\begin{definition}[Sections and quasi-sections]\ \par
Recall that $\pi_X: \Delta \to X $ is defined as $\pi_X(x,n)=x$. We say that $B_r(q) \times [-\pi/2, \pi/2]\subseteq \mathcal{M}$ is a section if $\pi_X: \pi^{-1} \big(B_r(q) \times [-\pi/2, \pi/2]\big) \to X$ is injective for any sufficiently small $r>0$. Further, $B_r(q) \times [-\pi/2, \pi/2] \subseteq \mathcal{M}$ is a quasi-section if for any sufficiently small $r>0$ there is a measurable set $S_r \subseteq B_r(q) \times [-\pi/2, \pi/2]$, such that  $\mu_{\mathcal{M}}[\big(B_r(q) \times [-\pi/2, \pi/2]\big)\setminus S_r]=o\big(\mu_{\mathcal{M}}\big[B_r(q) \times [-\pi/2, \pi/2]\big]\big)$,  and $\pi_X: \pi^{-1} S_r \to X$ is injective. In this case, we also refer to $S_r$ as a section in $B_r(q) \times [-\pi/2, \pi/2]$.  
\end{definition}
\begin{remark}\label{remarkonsection}
In applications to two-dimensional billiards, $B_r(q)\times [-\pi/2,\pi/2]$ is a strip in $\mathcal{M}$, and $ B_r(q)\times [-\pi/2,\pi/2]\setminus S_r$ is usually a union of finitely many rectangles, whose measures are of order $O(r^2)$. To avoid unnecessary complications, we always assume that $ B_r(q)\times [-\pi/2,\pi/2]\setminus S_r$ has a regular shape, e.g., as a union of finitely many rectangles.
\end{remark}

\begin{assumption}[\textbf{Geometric assumptions}]\label{assumption1}\ \par
\begin{enumerate}
    \item For a.e. $q \in \partial Q$ the set $B_r(q)\times [-\pi/2, \pi/2]$ is a quasi-section.
    \item $\bigcup_{i\ge 1}\partial X_i\subseteq \mathbb{S}$ (see the definition of $X_i$ in (\ref{parition})).
    \item There are constants $C>0$ and $\alpha \in (0,1]$, such that for any $\gamma^k$, where $k=u$ or $s$,  and for any $x,y \in \gamma^k$, \begin{align*}
        d_{f^j\gamma^k}(f^jx, f^jy)\le C d_{\gamma^k}(x,y)^\alpha \text{ for all }j\in[0,R(x)).
    \end{align*}
    (Observe that the condition $\bigcup_{i\ge 1}\partial X_i\subseteq \mathbb{S}$ implies that $\gamma^k \subseteq X_i$ for some $i\ge 1$).
    \item There exist two cones (i.e. affine subspaces) $C^u,C^s \subseteq \mathcal{T}(\mathcal{M})$, such that \[\dim(\interior{C^u}\bigcap \interior{C^s})<1, \quad (Df)C^u \subseteq C^u, \quad (Df)^{-1}C^s \subseteq C^s,\]
    and for all $n\ge 1$ and $\Leb_{\partial Q}$-a.e. $q \in \partial Q$
    \[\dim \big(\mathcal{T}(\{q\}\times [-\pi/2, \pi/2])\bigcap \interior{C^u}\big)<1, \quad \dim \big(\mathcal{T}(\{q\}\times [-\pi/2, \pi/2])\bigcap \interior{C^s}\big)<1,\]
    \[(Df)^n\mathcal{T}(\{q\}\times [-\pi/2, \pi/2]) \subseteq{\interior{C^u}}, \quad (Df)^{-n}\mathcal{T}(\{q\}\times [-\pi/2, \pi/2]) \subseteq{\interior{C^s}},\]where $\interior{C^u}$ (resp. $\interior{C^s}$) is the interior of $C^u$ (resp. $C^s$) and $\mathcal{T}(A)$ means a tangent space of a (sub)manifold $A$.
\end{enumerate}
\begin{remark}\ \par
\begin{enumerate}
    \item The existence of $C^u, C^s$ always holds for two dimensional hyperbolic billiards. 
    \item The H\"older condition is natural, and it is traditionally used for hyperbolic systems, (particularly, for billiards).
\end{enumerate}
 
\end{remark}
\end{assumption}

Now we can present the main results in this subsection.
\begin{theorem}[Poisson limit laws, see \cite{Subbb}]\label{thm}\ \par
Suppose that dynamical system $(\mathcal{M}, f)$ has a CMZ structure $(X, f^R)$ (see Definitions \ref{cmz} and \ref{inducesystem}), and the  Assumption \ref{assumption1} holds. Then, when $r\to 0$, we have $\mathcal{N}^{r,q} \to_d Poisson$ (see Definition \ref{pdef}) for $\Leb_{\partial Q}$-a.s. $q\in\partial Q$.
\end{theorem}
\begin{remark}
    The convergence rate (see Definition \ref{convergencerate}) for this approximation in Theorem \ref{thm} can be obtained by maximal large deviations \cite{mldp} and assuming polynomial mixing rates, while in Theorem \ref{thm} the convergence can be proven without specifying mixing rates.
\end{remark}

\begin{corollary}[First hitting times, see \cite{Subbb}]\ \par
Under the same conditions as in Theorem \ref{thm} consider the moment of time when the first hitting of (passage through) the hole occurs, i.e., $\tau_{r,q}(x):=\inf\left\{n \ge 1: f^n(x) \in B_r(q)\times [-\pi/2, \pi/2]\right\}$ for any $x \in \mathcal{M}$. Then for any $t>0$ and almost every $q \in \partial Q$, the following relation for the first hitting probability holds
\begin{equation*}
    \lim_{r\to 0}\mu_{\mathcal{M}}\Big\{\tau_{r,q}>t\big/\mu_{\mathcal{M}}\Big(B_r(q)\times [-\pi/2, \pi/2]\Big)\Big\}=e^{-t}.
\end{equation*}
\end{corollary}

We present now applications of Theorem \ref{thm} to some weakly chaotic billiards.
\subsubsection{Squashes or Stadium-type billiards}
\begin{figure}[!htb]
   \begin{minipage}{0.6\textwidth}
     \includegraphics[width=.6\linewidth]{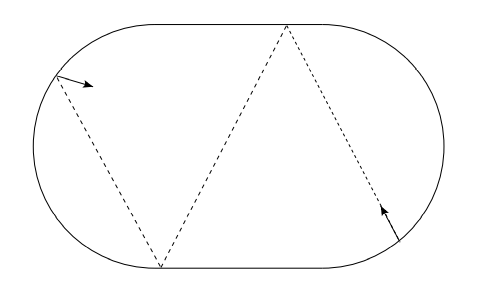}
     \caption{Stadium billiard}
     \label{F3}
   \end{minipage}\hfill
   \begin{minipage}{0.6\textwidth}
     \includegraphics[width=.6\linewidth]{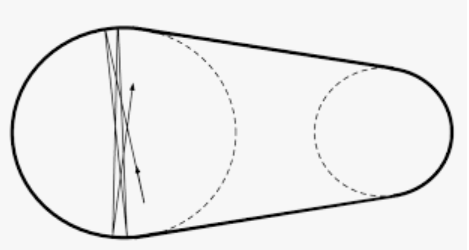}
     \caption{Squash billiard}\label{F4}
   \end{minipage}
\end{figure}

A billiard table $Q$ of a squash billiard is a convex domain bounded by two circular arcs and two straight (flat) segments tangent to the arcs at their common endpoints. A squash billiard is called a stadium if flat sides are parallel, see Figure \ref{F3}. Initially called squash billiards, they were later sometimes called ``skewed" stadia, drive-belt billiards, etc, see Figure \ref{F4}. Note that squashes contain a boundary arc, which is longer than a half circle.  Applying Theorem \ref{thm} we have the following.

\begin{corollary}
Theorem \ref{thm} holds for squash billiards and stadiums.
\end{corollary}

\subsubsection{Another class of billiards with focusing components} 
\begin{figure}[!htb]
   \begin{minipage}{0.6\textwidth}
\includegraphics[width=.6\linewidth]{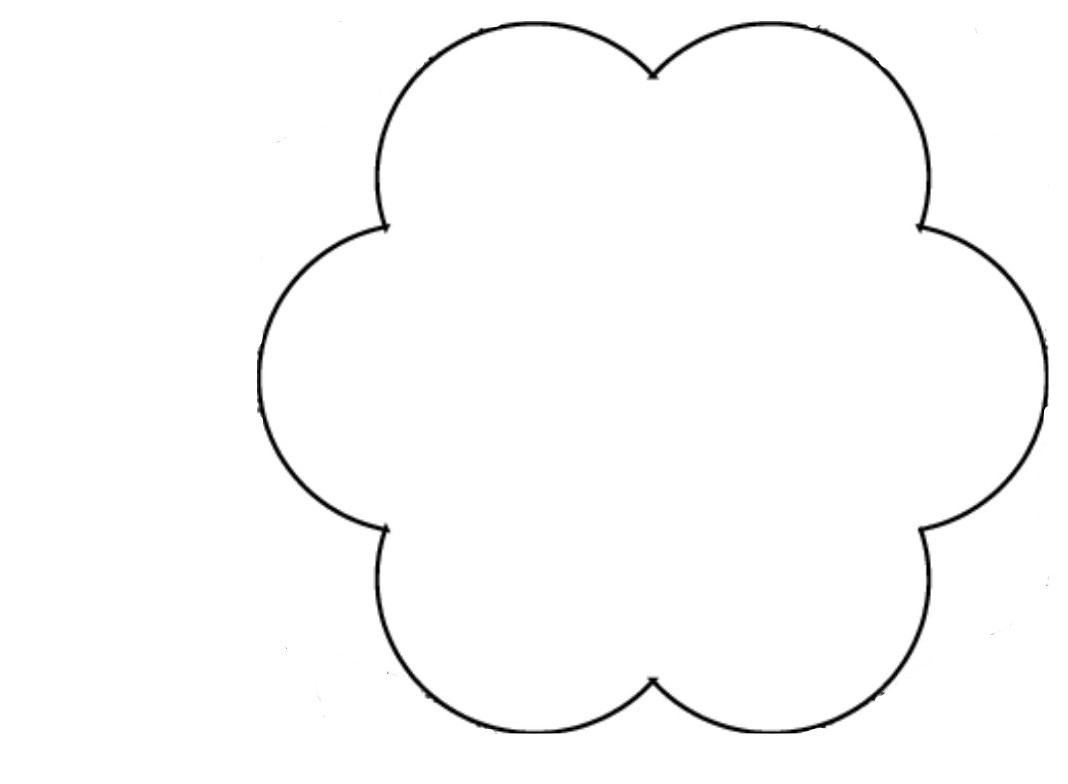}
     \caption{Flower billiard}
     \label{F5}
   \end{minipage}\hfill
   \begin{minipage}{0.48\textwidth}
     \includegraphics[width=.6\linewidth]{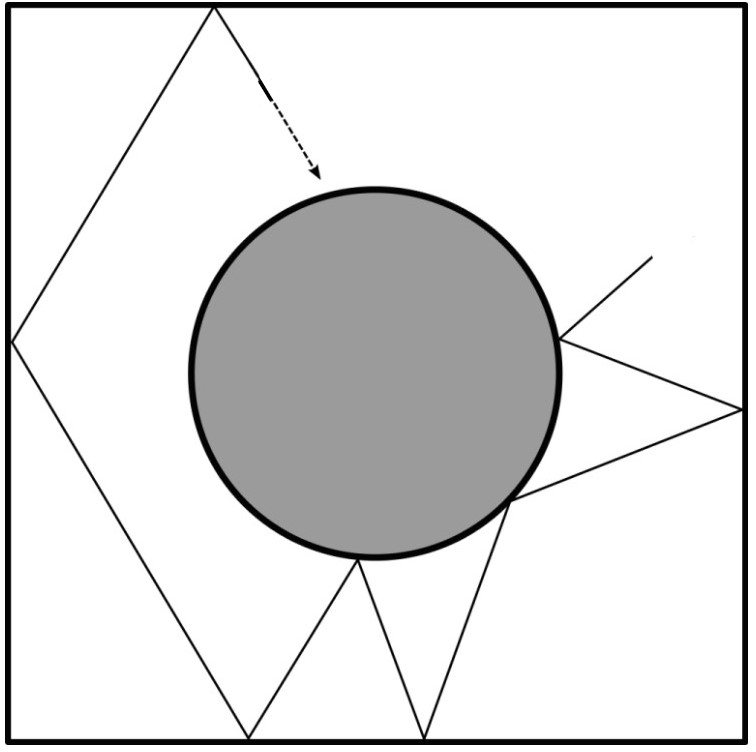}
     \caption{Semi-dispersing billiard}\label{F6}
   \end{minipage}
\end{figure}

In this section we consider billiard tables $Q$ for which each smooth component $\Gamma_i \subseteq \partial Q$ of the boundary is either dispersing, i.e., convex inwards, or focusing, i.e., convex
outwards. A curvature of every dispersing component is bounded
away from zero and infinity. We assume that every focusing component is an arc of a circle, and that there are no points of $\partial Q$ on that circle or inside it, other than the arc itself (the strong focusing or SFC-condition). We assume also that two dispersing components intersect (if they do) transversally (i.e., there are no cusps) and, besides, every focusing arc is not longer than a half of the corresponding circle, e.g. see the Figure \ref{F5}.
Applying Theorem \ref{thm} we have the following.

\begin{corollary}
Theorem \ref{thm} holds for the class of  billiards considered in this subsection.
\end{corollary}

\subsubsection{Semi-dispersing billiards} In this subsection we consider billiard tables of the following type. Let $R_0 \subseteq \mathbb{R}^2$ be a rectangle, and scatterers $B_1, \cdots, B_r \subseteq \interior R_0$ are open strictly convex sub-domains with smooth, (at least $C^3)$), or piece-wise smooth boundaries whose curvature is bounded away from zero,
and such that $B_i \bigcap B_j=\emptyset$ for $i\neq j$. The boundary of a billiard table $Q = R_0\setminus \bigcup_iB_i$ is partially dispersing (convex inwards) and partially neutral (flat), e.g. see Figure \ref{F6}. The flat part is $\partial R_0$. Applying Theorem \ref{thm} we have the following.
\begin{corollary}
Theorem \ref{thm} holds for the class of semi-dispersing billiards considered in this subsection.
\end{corollary}

\section{Central limit theorems and invariance principles}In this section we consider random dynamical systems.
A special subclass of random systems is formed by stationary systems.
\begin{definition}[Random and stationary dynamical systems]\ \par
     A random dynamical system is defined as $(\mathcal{M},(T_{\omega})_{\omega \in \Omega},  \Omega, \mathbb{P}, \sigma, (\mu_{\omega})_{\omega \in \Omega})$. Here  $\sigma: (\Omega, \mathbb{P})\to (\Omega,\mathbb{P})$ preserves a probability measure $\mathbb{P}$, and it is ergodic, $T_{\omega}: \mathcal{M} \to \mathcal{M}$ is a map,  $\mu_{\omega}$ is a probability measure on $X$ satisfying $(T^k_{\omega})_{*}\mu_{\omega}=\mu_{\sigma^k \omega}, k \ge 0$ where $T_{\omega}^k:=T_{\sigma^{k-1}\omega}\circ \cdots \circ T_{\sigma \omega}\circ T_{\omega}$.

     If $\Omega$ is a single point, then we call this system a stationary dynamical system.
\end{definition}

\begin{definition}[Vector-valued Almost Sure Invariance Principle (VASIP)]\
\par\label{rds}For a random dynamical system $(\mathcal{M},(T_{\omega})_{\omega \in \Omega},  \Omega, \mathbb{P}, \sigma, (\mu_{\omega})_{\omega \in \Omega})$, consider observables (functions on the phase space of all states of the system) $ \{\phi_{\omega} \in L^{\infty}(\mathcal{M},\mu_\mathcal{M}; \mathbb{R}^d): \omega \in \Omega\}$ which satisfy $$\sup_{\omega \in \Omega}||\phi_{\omega}||_{\infty}< \infty ,\quad  \int{\phi_{\omega}  d\mu_{\omega}}=0 \text{ for all } \omega \in \Omega.$$

We say that $ (\phi_{\sigma^k\omega} \circ T_{\omega}^k)_{k \ge 1, \omega \in \Omega} $ satisfies a quenched vector-valued almost sure invariance principle (quenched VASIP) if, for a.e. $\omega \in \Omega$, there exists a constant $\epsilon \in (0,1) $ and a family of zero mean $d$-dimensional independent Gaussian random vectors $ (G^{\omega}_k)_{k \ge 1}$ defined on some extended probability space of $(\mathcal{M}, \mu_{\omega})$ such that
\begin{equation} \label{qmatching}
\sum_{k=1}^{n} \phi_{\sigma^k\omega} \circ T_{\omega}^k - \sum_{k=1}^{n} G^{\omega}_k=o(n^{\frac{1-\epsilon}{2}}) \text{ a.s.}
\end{equation}
\begin{equation} \label{qvariancegrowth}
 \sigma^2_n(\omega)=\sum_{k=1}^{n} \tilde{\mathbb{E}}^{\omega}[{ (G^{\omega}_k) (G^{\omega}_k)^{T}}]+ o(n^{1-\epsilon}),
\end{equation}
\begin{equation}\label{qlinear}
\sigma_n^2(\omega) \approx_{\omega} n  I_{d \times d},
\end{equation}
where $\sigma_n^2(\omega):=\int(\sum_{k=1}^{n} \phi_{\sigma^k\omega} \circ T_{\omega}^k )(\sum_{k=1}^{n} \phi_{\sigma^k \omega}\circ T_{\omega}^k )^{T} d\mu_{\omega}$ and $\tilde{\mathbb{E}}^{\omega}$ in (\ref{qvariancegrowth}) is the expectation w.r.t. the probability $\tilde{P}^{\omega}$ on an extended probability space of $(\mathcal{M},  \mu_{\omega})$. Here the extended probability space is $(\mathcal{M}, \mu_{\omega}) \otimes \otimes_{ i \ge 1}(I_i, \Leb) $ where $I_i$ is a copy of $[0,1]$.
\end{definition}

\begin{remark} If $d=1$, then the quenched VASIP implies the quenched central limit theorem (CLT): \begin{align*}
    \lim_{n\to \infty }\mu_{\omega} \Big(\frac{\sum_{k=1}^{n} \phi_{\sigma^k \omega} \circ T_{\omega}^k}{\Sigma \sqrt{n} }\le z\Big) =\frac{1}{\sqrt{2\pi}}\int_{-\infty}^{z}e^{-s^2/2}ds \text{ for all } z \in \mathbb{R}, 
\end{align*}(in other words, $\frac{\sum_{k=1}^{n} \phi_{\sigma^k \omega} \circ T_{\omega}^k}{\Sigma \sqrt{n} }$ converges in distribution to a Gaussian variable $\mathcal{N}(0,1)$) and the quenched iterated logarithm  law(LIL): \begin{align*}
  \limsup_{n\to \infty}\frac{\sum_{k=1}^{ n} \phi_{\sigma^k\omega} \circ T_{\omega}^k}{\sqrt{n\log\log n}}=\Sigma \quad \mu_{\omega}\text{-a.s.}
\end{align*}
for a.e. $\omega \in \Omega$ and a constant $\Sigma > 0$, which does not depend on $\omega$. 
\end{remark}

\begin{remark}\label{singlept}
    If $\Omega$ is a single point $\{\omega_0\}$, then the random system becomes a stationary system. Then VASIP, CLT and LIL have the following forms, respectively: 
\begin{gather*}
    \sum_{k=1}^{n} \phi_{\omega_0} \circ T_{\omega_0}^k - \sum_{k=1}^{n} G_k=o(n^{\frac{1-\epsilon}{2}}) \text{ a.s.}\\
     \lim_{n\to \infty }\mu_{\omega_0} \Big(\frac{\sum_{k=1}^{n} \phi_{\omega_0} \circ T_{\omega_0}^k}{\Sigma \sqrt{n} }\le z\Big) =\frac{1}{\sqrt{2\pi}}\int_{-\infty}^{z}e^{-s^2/2}ds \text{ for all } z \in \mathbb{R}, \\
     \limsup_{n\to \infty}\frac{\sum_{k=1}^{ n} \phi_{\omega_0} \circ T_{\omega_0}^k}{\sqrt{n\log\log n}}=\Sigma \quad \mu_{\omega_0}\text{-a.s.}
\end{gather*}
\end{remark}

The following is a general VASIP result for random dynamical systems.
\begin{theorem}[Random Dynamical Systems, see \cite{Sutams}]\ \label{thm3} \par
Consider a random dynamical system $(\mathcal{M}, (T_{\omega})_{\omega \in \Omega},  \Omega, \mathbb{P}, \sigma, (\mu_{\omega})_{\omega \in \Omega})$, where $d\mu_{\omega}=h_{\omega}d\mu_\mathcal{M}$ is a probability distribution  on $\mathcal{M}$ and $\mu_{\mathcal{M}}$ is a measure on $\mathcal{M}$. Let the transfer operators (Perron-Frobenius operators)  $P_{\omega}$ associated to $T_{\omega}$ are defined by the  following duality relations
\begin{equation*}
    \int g \cdot P_{\omega} (h) d\mu_{\mathcal{M}}= \int g \circ T_{\omega} \cdot h d\mu_{\mathcal{M}} \text{ for all } h \in L^1(\mu_{\mathcal{M}})\text{, } g \in L^{\infty}(\mu_{\mathcal{M}}).
\end{equation*}

Let the observables $\{\phi_{\omega} \in L^{\infty}(\mathcal{M},\mu_{\mathcal{M}};\mathbb{R}^d):  \omega \in \Omega\}$ be the ones from the Definition \ref{rds}. Suppose that they satisfy the conditions (\ref{A1'})-(\ref{A3'}) below, and there is a constant $C>0$ such that for any $i,j\ge 0,n \ge 1$, and for almost every $\omega\in \Omega$,
\begin{equation} \label{A1'}
\int |P_{\sigma^{i}\omega}^{n} (\phi_{\sigma^i \omega}  h_{\sigma^{i} \omega})|d\mu_{\mathcal{M}} \le C n^{1-1/\alpha}, \tag{A1'}
\end{equation}
\begin{equation} \label{A2'}
\int \Big|P_{\sigma^{i}\omega}^{n} \Big[(\phi_{\sigma^i \omega}  \phi_{\sigma^i \omega}^T - \int \phi_{\sigma^i \omega}   \phi_{\sigma^i \omega}^T  d\mu_{\sigma^i\omega}) h_{\sigma^{i}\omega}\Big] \Big|  d\mu_{\mathcal{M}} \le C  n^{1-1/\alpha}, \tag{A2'}
\end{equation}
\begin{equation} \label{A3'}
\int \Big|P_{\sigma^{i+j} \omega}^{n} \Big[P^{j}_{\sigma^{i} \omega}(\phi_{\sigma^{i} \omega} h_{\sigma^{i}\omega})  \phi_{\sigma^{i+j} \omega}^T - h_{\sigma^{i+j}\omega}\int P^{j}_{\sigma^{i} \omega}(\phi_{\sigma^{i} \omega}  h_{\sigma^{i}\omega}) \phi_{\sigma^{i+j} \omega}^T d\mu_{\mathcal{M}}\Big] \Big| d\mu_{\mathcal{M}}
\le C n^{1-1/\alpha}, \tag{A3'}
\end{equation}
where $\alpha < 1/2$, $|\cdot|$ is the Euclidean norm for vectors or the spectral norm for matrices.

Then there are two linear subspaces (which do not depend on $\omega$): $W_1, W_2 \subseteq{\mathbb{R}^d}$, $\mathbb{R}^d= W_1 \bigoplus W_2$ with projections $\pi_1:W_1\bigoplus W_2 \to W_1, \pi_2: W_1 \bigoplus W_2 \to W_2$ such that
\begin{enumerate}
    \item $(\pi_1 \circ  \phi_{\sigma^k \omega}  \circ T^{k}_{\omega})_{k \ge 1, \omega \in \Omega}$ satisfies the quenched VASIP (see Definition \ref{rds}).
    \item $(\pi_2 \circ  \phi_{\sigma^k \omega}  \circ T^{k}_{\omega})_{k \ge 1, \omega \in \Omega}$ is a coboundary, i.e., there is $\psi \in L^1(\Omega \times X, d\mu_{\omega}d\mathbb{P})$ such that \[  \pi_2 \circ \phi_{\sigma \omega}(T_{\omega}x)=\psi(\sigma \omega,T_{\omega}x)-\psi(\omega,x) \text{ a.e. }(\omega,x)\in \Omega \times X.\]
\end{enumerate}
\end{theorem}

The paper \cite{LSV} considers the so called Pomeau-Manneville type maps, which are expanding everywhere besides at an indifferent (parabolic) point. Existence of just one such point may essentially change  dynamics. This class of maps is defined in the following way. Let $\beta\in [0,1)$,

\begin{eqnarray}\label{pmmap}
 T_{\beta}(x) =
\begin{cases}
x+2^{\beta}x^{1+\beta},      & 0\le x \le 1/2\\
2x-1,  & 1/2< x \le 1 \\
\end{cases}.
\end{eqnarray}

According to \cite{LSV, Young2}, its mixing rate is $O(n^{1-1/\beta})$ and this map preserves a probability measure absolutely continuous respect to the Lebesgue measure on $[0,1]$. If we randomly perturbed the parameters $\beta \in [0,1/2)$, and apply Theorem \ref{thm3} to the class of these random systems, then the following result holds.

\begin{corollary}[Polynomially mixing random systems, see \cite{Sutams}]\ \label{cor7} \par
Let $\mu_{\mathcal{M}}:=\Leb$ be the Lebesgue measure on $\mathcal{M}:=[0,1]$,  $\Omega:=[0,1/2)^{\mathbb{Z}}$, and $T_{\omega}:=T_{\omega_{[0]}}$ are the Pomeau-Manneville type maps (\ref{pmmap}) which are picked from $\{T_{\beta}: \beta \in [0,1/2)\}$. Let also $ \sigma : \Omega \to \Omega$ is a $\mathbb{P}$-preserving ergodic and invertible left shift, and  $\omega_{[0]}$ is the 0th entry of $\omega$.

Then there are functions $h_{\omega} \in L^1(\Leb)$ and probabilities $d\mu_{\omega}: =h_{\omega}d\Leb$ such that $(T_{\omega})_{*} \mu_{\omega}=\mu_{\sigma \omega}$ for a.e. $\omega \in \Omega$. Moreover, suppose that the observables $(\phi_{\omega})_{\omega \in \Omega} \subseteq \mathrm{Lip}([0,1]; \mathbb{R}^d)$ satisfy the relations $\sup_{\omega}||\phi_{\omega}||_{\mathrm{Lip}}<\infty$ and $\int \phi_{\omega}d\mu_{\omega}=0$. Then there are two linear subspaces (which do not depend on $\omega$) $W_1, W_2 \subseteq{\mathbb{R}^d}$, $\mathbb{R}^d= W_1 \bigoplus W_2$ with projections $\pi_1:W_1\bigoplus W_2 \to W_1, \pi_2: W_1 \bigoplus W_2 \to W_2$ such that 

\begin{enumerate}
    \item $(\pi_1 \circ  \phi_{\sigma^k \omega}  \circ T^{k}_{\omega})_{k \ge 1,\omega \in \Omega}$ satisfies the quenched VASIP.
    \item $(\pi_2 \circ  \phi_{\sigma^k \omega}  \circ T^{k}_{\omega})_{k \ge 1,\omega \in \Omega}$ is a coboundary, which means that there is $\psi \in L^1(\Omega \times [0,1], d\mu_{\omega}d\mathbb{P})$ such that \[\pi_2 \circ  \phi_{\sigma \omega}(T_{\omega}x)=\psi(\sigma \omega,T_{\omega}x)-\psi(\omega,x) \text{ a.e. }(\omega,x)\in \Omega \times [0,1].\]
\end{enumerate}

\end{corollary}
\begin{remark}
$\sigma: \Omega \to \Omega$ is ergodic, and $1/2$ is optimal, see also \cite{nonclt}. If $\sigma$ is Bernoulli, then the range of $\beta$ can be larger than $[0,1/2)$ (see  \cite{Suetds}).
\end{remark}
\begin{remark}\label{singlemaplimit}
By Remark \ref{singlept}, when $\beta\in [0,1/2), \Omega=\{\beta\}$ and $\phi$ has zero mean, then  we obtain the corresponding results for the  single dynamics $T_{\beta}$ and a zero mean function $\phi$  \begin{gather*}
    \sum_{k=1}^{n} \phi \circ T_{\beta}^k - \sum_{k=1}^{n} G_k=o(n^{\frac{1-\epsilon}{2}}) \text{ a.s.}\\
    \frac{\sum_{k=1}^{n} \phi \circ T_{\beta}^k}{\Sigma \sqrt{n} } \to_d \mathcal{N}(0,1)
\end{gather*} (Gaussian distribution) or there is a function $\chi$ such that $\phi=\chi\circ T_{\beta}-\chi$ a.s.
\end{remark}

If $\beta \in (1/2,1)$, then weakly chaotic stationary dynamics $T_{\beta}$ does not satisfy the standard CLT or VASIP,  but the following different limit law holds \cite{nonclt}.

\section{Stable limit laws}
Let $1/2<\beta <1$, and the corresponding stationary dynamics $T_{\beta} : [0, 1] \to [0, 1]$ be the Liverani-Saussol-Vaienti map described by (\ref{pmmap}). 

\begin{theorem}[see \cite{nonclt}]\label{gouezelstable} \par
   Let $h$ be the density of its unique absolutely continuous invariant probability measure $\mu_{\beta}$, and $\phi: [0, 1] \to \mathbb{R}$ is H\"older with  $\int \phi d\mu_{\beta} = 0$.
 
If $\phi (0) \neq 0$, then the following stable limit law holds: \[\frac{\sum_{i \le n} \phi\circ T_{\beta}^i}{n^\beta} \to_d X_{1/\beta,c,sgn(\phi(0))}\] whose characteristic function is $\mathbb{E}(e^{itX_{1/\beta,c,sgn(\phi(0))}} ) = e^{-c|t|^{1/\beta}(1-isgn(\phi(0))sgn(t)\tan(\pi/(2\beta) ))} $ with $c = \frac{h(1/2)}{4(\beta/|\phi (0)|)^{1/\beta}} \Gamma(1-1/\beta)\cos(\pi/(2\beta))$.

If $\phi(0) = 0$, we may also assume that $|\phi(x)|\le Cx^{\gamma}$ with $\gamma >\beta-1/2$. Then there
exists a constant $\sigma^2\ge 0$ such that \[\frac{\sum_{i\le n}\phi\circ T_{\beta}^i}{\sqrt{n}} \to_d \mathcal{N}(0, \sigma^2).\] Moreover,
$\sigma^2 = 0$ if and only if there exists a measurable function $\chi$ such that $\phi = \chi \circ T_{\beta} -\chi$.
\end{theorem}
Besides expadning maps $T_{\beta}$,  we next consider two-dimensional dispersing billiards with a cusp at a flat point. More precisely, for any fixed constant $\eta> 2$
(which will determine the sharpness of the cusp), we take a dispersing
billiard table $Q = Q_{\eta}$ with a boundary consisting of a finite number of $C^3$
smooth curves $\Gamma_i$, $i = 1 \cdots, n_0$, where $n_0 \ge 3$, where a cusp formed by two  curves $\Gamma_1, \Gamma_2$ such that  cusp at $P = \Gamma_1 \bigcap \Gamma_2$. Moreover, both $\Gamma_1$ and $\Gamma_2$ have zero derivatives up to $\eta- 1$ order
at $P$, and the $\eta$-order derivative is not zero. We also assume that all other
boundary components are dispersing and have curvature bounded away from
zero.
Choose now a Cartesian coordinate system $(s, z)$ with the origin at $P$, and with
the horizontal $s$-axis being the tangent line to both $\Gamma_1$ and $\Gamma_2$. Assume also that $\Gamma_1$ and
$\Gamma_2$ can be represented as
\begin{equation}\label{cusps}
    z_1(s) = \eta^{-1}s^\eta, z_2(s) = -\eta^{-1}s^\eta
\end{equation}
for $s \in [0, \epsilon_0]$ with $\epsilon_0 > 0$ being a small fixed number. The mixing rate for this billiard system is $O(n^{-1/(\eta-1)})$. A stable limit law for this billiard system ($\mathcal{M}, f$) (see Definition \ref{billiarddef}) is proved by Paul and Zhang \cite{hongkunstable}.

\begin{theorem}[see \cite{hongkunstable}]\ \par
    Let $Q_{\eta}$, where
$\eta \in(2,\infty)$, is a billiard table with a cusp defined by (\ref{cusps}), $\alpha = \eta/(\eta-1)$ and suppose that a mean zero function $\psi$ is  H\"older and defined on  $\mathcal{M}$. Let also $$I_\psi:= 1/4 \int^{\pi}_0(\psi(q''
, \phi) + \psi(q', \phi)) \sin^{1/\alpha} \phi d\phi \neq 0,$$ where $q', q''$ are the $q$-coordinates of the cusps point $P$ in the phase space $\mathcal{M}$ with respect to $\Gamma_1$ and $\Gamma_2$, respectively. Then as $n \to \infty$, \[\frac{\sum_{i \le n}\psi \circ f^i}{n^{1/\alpha}} \to_d S_{\alpha, \sigma}\]
where $\sigma = \frac{2I_\psi}{\eta \Leb_{\partial Q}\partial Q}$ and the characteristic function of $S_{\alpha,\sigma}$ is $\mathbb{E}(e^{itS_{\alpha, \sigma}})=e^{-|\sigma|^{\alpha}|t|^{\alpha}(1-isgn(t)\tan(\alpha\pi/2 ))} $.
\end{theorem}

\begin{remark}
    We note that when the mixing rate is slower and slower, the Poisson limit laws (Theorem \ref{thm}) still hold, while the central limit theorems (Remark \ref{singlemaplimit})  will change to stable limit laws (Theorem \ref{gouezelstable}). 
\end{remark}

\begin{figure}[!h]
    \centering
    \includegraphics[width=90mm]{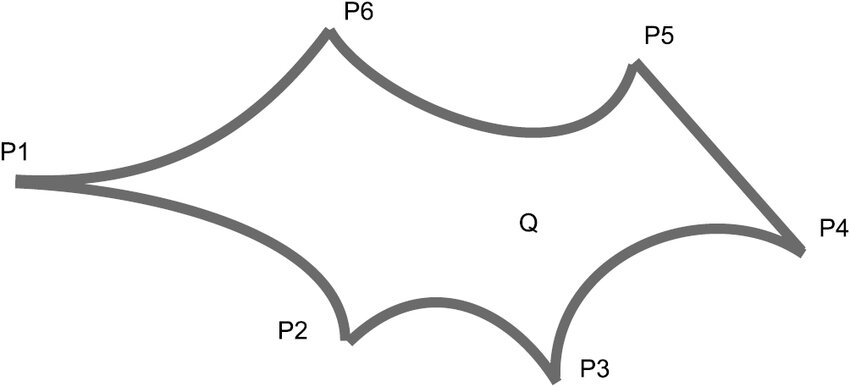}
    \caption{A billiard table with several cusps at flat-points \cite{cuspsmorethanone}}\label{199}
\end{figure}
\begin{remark}
       In a slightly more complicated case, a billiard table has, possibly non-isometric and asymmetric cusps, at flat points  (see Figure \ref{199}). It was shown in\cite{cuspsmorethanone, melbournelevy1,melbournelevy2} that in such situation a functional limit law with H\"older observables holds in Skorokhod $\mathcal{M}_1$-topology (or the weaker $\mathcal{M}_2$-topology) if $$I_\psi(s):= 1/4 \int^{s}_0(\psi(q''
, \phi) + \psi(q', \phi)) \sin^{1/\alpha} \phi d\phi, s\in [0,\pi]$$ is monotone (or the image of $I_\psi$ is contained in $ [0, I_\psi(\pi)]$). We refer the readers to \cite{cuspsmorethanone, melbournelevy1,melbournelevy2} for more details.
\end{remark}

\begin{acknowledgements}
L.B. was partially supported by the NSF grant DMS-2054659. Y. Su thanks Prof. Ian Melbourne for helpful comments.

\end{acknowledgements}


\medskip

% Non-BibTeX users  please use
\bibliography{bibtext}

@article {melbournelevy2,
    AUTHOR = {Jung, Paul and Melbourne, Ian and P\`ene, Fran\c coise and
              Varandas, Paulo and Zhang, Hong-Kun},
     TITLE = {Necessary and sufficient condition for M2-convergence to a {L}\'evy process for billiards with
              cusps at flat points},
   JOURNAL = {Stoch. Dyn.},
  FJOURNAL = {Stochastics and Dynamics},
    VOLUME = {21},
      YEAR = {2021},
    NUMBER = {5},
     PAGES = {Paper No. 2150024, 8},
      ISSN = {0219-4937,1793-6799},
   MRCLASS = {37C83 (37A50 37D25 60F17)},
  MRNUMBER = {4297947},
       DOI = {10.1142/S0219493721500246},
       URL = {https://doi.org/10.1142/S0219493721500246},
}

@article {melbournelevy1,
    AUTHOR = {Melbourne, Ian and Varandas, Paulo},
     TITLE = {Convergence to a {L}\'evy process in the Skorohod
             M1 and M2 topologies for nonuniformly
              hyperbolic systems, including billiards with cusps},
   JOURNAL = {Comm. Math. Phys.},
  FJOURNAL = {Communications in Mathematical Physics},
    VOLUME = {375},
      YEAR = {2020},
    NUMBER = {1},
     PAGES = {653--678},
      ISSN = {0010-3616,1432-0916},
   MRCLASS = {37C83 (37A50 37D25 60F17)},
  MRNUMBER = {4082176},
MRREVIEWER = {Hongkun\ Zhang},
       DOI = {10.1007/s00220-019-03501-9},
       URL = {https://doi.org/10.1007/s00220-019-03501-9},
}

@article {mldp,
    AUTHOR = {Bunimovich, Leonid A. and Su, Yaofeng},
     TITLE = {Maximal large deviations and slow recurrences in weakly
              chaotic systems},
   JOURNAL = {Adv. Math.},
  FJOURNAL = {Advances in Mathematics},
    VOLUME = {432},
      YEAR = {2023},
     PAGES = {Paper No. 109267, 58},
      ISSN = {0001-8708,1090-2082},
   MRCLASS = {37A50 (37A25)},
  MRNUMBER = {4630843},
MRREVIEWER = {Ivan\ Podvigin},
       DOI = {10.1016/j.aim.2023.109267},
       URL = {https://doi.org/10.1016/j.aim.2023.109267},
}

@article {K,
    AUTHOR = {Kolmogorov, Andrei Nikolaevich},
     TITLE = {A new metric invariant of transient dynamical systems and automorphisms in Lebesgue spaces},
   JOURNAL = {Dokl. Akad. Nauk SSSR},
  FJOURNAL = {Doklady Akademii Nauk SSSR},
    VOLUME = {119},
      YEAR = {1958},
     PAGES = {861--864},
      ISSN = {5},
       URL = {http://mi.mathnet.ru/dan22922},
}

@article {vaientiadv,
    AUTHOR = {Alves, Jos\'{e} F. and Freitas, Jorge M. and Luzzatto, Stefano and
              Vaienti, Sandro},
     TITLE = {From rates of mixing to recurrence times via large deviations},
   JOURNAL = {Adv. Math.},
  FJOURNAL = {Advances in Mathematics},
    VOLUME = {228},
      YEAR = {2011},
    NUMBER = {2},
     PAGES = {1203--1236},
      ISSN = {0001-8708},
   MRCLASS = {37C40 (37A25 37D25)},
  MRNUMBER = {2822221},
MRREVIEWER = {Stefano Galatolo},
       DOI = {10.1016/j.aim.2011.06.014},
       URL = {https://doi.org/10.1016/j.aim.2011.06.014},
}

@article {hongkunstable,
    AUTHOR = {Jung, Paul and Zhang, Hong-Kun},
     TITLE = {Stable laws for chaotic billiards with cusps at flat points},
   JOURNAL = {Ann. Henri Poincar\'{e}},
  FJOURNAL = {Annales Henri Poincar\'{e}. A Journal of Theoretical and
              Mathematical Physics},
    VOLUME = {19},
      YEAR = {2018},
    NUMBER = {12},
     PAGES = {3815--3853},
      ISSN = {1424-0637,1424-0661},
   MRCLASS = {37C83},
  MRNUMBER = {3877517},
       DOI = {10.1007/s00023-018-0726-y},
       URL = {https://doi.org/10.1007/s00023-018-0726-y},
}

@article {Kolmo,
    AUTHOR = {Kolmogorov, A. N.},
     TITLE = {A new metric invariant of transient dynamical systems and
              automorphisms in {L}ebesgue spaces},
   JOURNAL = {Dokl. Akad. Nauk SSSR (N.S.)},
  FJOURNAL = {Dokl. Akad. Nauk SSSR (N.S.)},
    VOLUME = {119},
      YEAR = {1958},
     PAGES = {861--864},
   MRCLASS = {28.00 (94.00)},
  MRNUMBER = {103254},
}

@article {Sutams,
    AUTHOR = {Su, Yaofeng},
     TITLE = {Vector-valued almost sure invariance principles for
              (non)stationary and random dynamical systems},
   JOURNAL = {Trans. Amer. Math. Soc.},
  FJOURNAL = {Transactions of the American Mathematical Society},
    VOLUME = {375},
      YEAR = {2022},
    NUMBER = {7},
     PAGES = {4809--4848},
      ISSN = {0002-9947,1088-6850},
   MRCLASS = {37C99 (60F17)},
  MRNUMBER = {4439492},
MRREVIEWER = {Jianyu\ Chen},
       DOI = {10.1090/tran/8609},
       URL = {https://doi.org/10.1090/tran/8609},
}

@article {LSV,
    AUTHOR = {Liverani, Carlangelo and Saussol, Beno\^{\i}t and Vaienti, Sandro},
     TITLE = {A probabilistic approach to intermittency},
   JOURNAL = {Ergodic Theory Dynam. Systems},
  FJOURNAL = {Ergodic Theory and Dynamical Systems},
    VOLUME = {19},
      YEAR = {1999},
    NUMBER = {3},
     PAGES = {671--685},
      ISSN = {0143-3857},
   MRCLASS = {37D25 (28D05 37A99 37C40 37E05)},
  MRNUMBER = {1695915},
MRREVIEWER = {Makoto Mori},
       DOI = {10.1017/S0143385799133856},
       URL = {https://doi.org/10.1017/S0143385799133856},
}

@ARTICLE{explicit,
    AUTHOR = {Bunimovich, Leonid A. and Su, Yaofeng},
     TITLE = {Estimates of constants in the limit theorems for chaotic
              dynamical systems},
   JOURNAL = {Stoch. Dyn.},
  FJOURNAL = {Stochastics and Dynamics},
    VOLUME = {24},
      YEAR = {2024},
    NUMBER = {1},
     PAGES = {Paper No. 2450004, 15},
      ISSN = {0219-4937,1793-6799},
   MRCLASS = {37A50 (37D45)},
  MRNUMBER = {4745238},
       DOI = {10.1142/S0219493724500047},
       URL = {https://doi.org/10.1142/S0219493724500047},
}

@article {melbourne09,
    AUTHOR = {Melbourne, Ian},
     TITLE = {Large and moderate deviations for slowly mixing dynamical
              systems},
   JOURNAL = {Proc. Amer. Math. Soc.},
  FJOURNAL = {Proceedings of the American Mathematical Society},
    VOLUME = {137},
      YEAR = {2009},
    NUMBER = {5},
     PAGES = {1735--1741},
      ISSN = {0002-9939},
   MRCLASS = {37D25 (28D05 37A05 37D50 60F10)},
  MRNUMBER = {2470832},
MRREVIEWER = {Jean-Ren\'{e} Chazottes},
       DOI = {10.1090/S0002-9939-08-09751-7},
       URL = {https://doi.org/10.1090/S0002-9939-08-09751-7},
}

@article {haydncmp,
    AUTHOR = {Haydn, Nicolai and Vaienti, Sandro},
     TITLE = {Limiting entry and return times distribution for arbitrary
              null sets},
   JOURNAL = {Comm. Math. Phys.},
  FJOURNAL = {Communications in Mathematical Physics},
    VOLUME = {378},
      YEAR = {2020},
    NUMBER = {1},
     PAGES = {149--184},
      ISSN = {0010-3616},
   MRCLASS = {37A50 (37D10 60B10 60E05)},
  MRNUMBER = {4124984},
       DOI = {10.1007/s00220-020-03795-0},
       URL = {https://doi.org/10.1007/s00220-020-03795-0},
}

@article {Subbb,
    AUTHOR = {Bunimovich, Leonid A. and Su, Yaofeng},
     TITLE = {Back to boundaries in billiards},
   JOURNAL = {Comm. Math. Phys.},
  FJOURNAL = {Communications in Mathematical Physics},
    VOLUME = {405},
      YEAR = {2024},
    NUMBER = {6},
     PAGES = {Paper No. 140, 74},
      ISSN = {0010-3616,1432-0916},
   MRCLASS = {37C83 (37A50 60F05)},
  MRNUMBER = {4751722},
       DOI = {10.1007/s00220-024-05002-w},
       URL = {https://doi.org/10.1007/s00220-024-05002-w},
}

@article {Chernovjsp,
    AUTHOR = {Chernov, N.},
     TITLE = {Decay of correlations and dispersing billiards},
   JOURNAL = {J. Statist. Phys.},
  FJOURNAL = {Journal of Statistical Physics},
    VOLUME = {94},
      YEAR = {1999},
    NUMBER = {3-4},
     PAGES = {513--556},
      ISSN = {0022-4715},
   MRCLASS = {37D50 (28D05 37A99 37C40 82C40)},
  MRNUMBER = {1675363},
MRREVIEWER = {Guido Gentile},
       DOI = {10.1023/A:1004581304939},
       URL = {https://doi.org/10.1023/A:1004581304939},
}

@article {CZnon,
    AUTHOR = {Chernov, N. and Zhang, H.-K.},
     TITLE = {Billiards with polynomial mixing rates},
   JOURNAL = {Nonlinearity},
  FJOURNAL = {Nonlinearity},
    VOLUME = {18},
      YEAR = {2005},
    NUMBER = {4},
     PAGES = {1527--1553},
      ISSN = {0951-7715},
   MRCLASS = {37D50 (37A25 37N20 82B05)},
  MRNUMBER = {2150341},
MRREVIEWER = {Roberto Markarian},
       DOI = {10.1088/0951-7715/18/4/006},
       URL = {https://doi.org/10.1088/0951-7715/18/4/006},
}

@article {collet,
    AUTHOR = {Chazottes, J.-R. and Collet, P.},
     TITLE = {Poisson approximation for the number of visits to balls in
              non-uniformly hyperbolic dynamical systems},
   JOURNAL = {Ergodic Theory Dynam. Systems},
  FJOURNAL = {Ergodic Theory and Dynamical Systems},
    VOLUME = {33},
      YEAR = {2013},
    NUMBER = {1},
     PAGES = {49--80},
      ISSN = {0143-3857},
   MRCLASS = {37D50 (60G10 62E17)},
  MRNUMBER = {3009103},
MRREVIEWER = {Tam\'{a}s I. Wiandt},
       DOI = {10.1017/S0143385711000897},
       URL = {https://doi.org/10.1017/S0143385711000897},
}

@article {Suetds,
    AUTHOR = {Su, Yaofeng},
     TITLE = {Random {Y}oung towers and quenched limit laws},
   JOURNAL = {Ergodic Theory Dynam. Systems},
  FJOURNAL = {Ergodic Theory and Dynamical Systems},
    VOLUME = {43},
      YEAR = {2023},
    NUMBER = {3},
     PAGES = {971--1003},
      ISSN = {0143-3857,1469-4417},
   MRCLASS = {37H30 (60F17)},
  MRNUMBER = {4544151},
       DOI = {10.1017/etds.2021.164},
       URL = {https://doi.org/10.1017/etds.2021.164},
}

@article {Sucmp,
    AUTHOR = {Bunimovich, Leonid A. and Su, Yaofeng},
     TITLE = {Poisson approximations and convergence rates for hyperbolic
              dynamical systems},
   JOURNAL = {Comm. Math. Phys.},
  FJOURNAL = {Communications in Mathematical Physics},
    VOLUME = {390},
      YEAR = {2022},
    NUMBER = {1},
     PAGES = {113--168},
      ISSN = {0010-3616},
   MRCLASS = {37D20 (37A50 60F05 60G55)},
  MRNUMBER = {4381186},
MRREVIEWER = {Ivan Podvigin},
       DOI = {10.1007/s00220-022-04309-w},
       URL = {https://doi.org/10.1007/s00220-022-04309-w},
}

@article {Alves,
    AUTHOR = {Alves, Jos\'{e} F. and Azevedo, Davide},
     TITLE = {Statistical properties of diffeomorphisms with weak invariant
              manifolds},
   JOURNAL = {Discrete Contin. Dyn. Syst.},
  FJOURNAL = {Discrete and Continuous Dynamical Systems. Series A},
    VOLUME = {36},
      YEAR = {2016},
    NUMBER = {1},
     PAGES = {1--41},
      ISSN = {1078-0947},
   MRCLASS = {37A05 (37C40 37D25 60F10)},
  MRNUMBER = {3369212},
MRREVIEWER = {Ricardo G\'{o}mez},
       DOI = {10.3934/dcds.2016.36.1},
              URL = {https://doi.org/10.3934/dcds.2016.36.1},
}

@article {Young1,
    AUTHOR = {Young, Lai-Sang},
     TITLE = {Statistical properties of dynamical systems with some
              hyperbolicity},
   JOURNAL = {Ann. of Math. (2)},
  FJOURNAL = {Annals of Mathematics. Second Series},
    VOLUME = {147},
      YEAR = {1998},
    NUMBER = {3},
     PAGES = {585--650},
      ISSN = {0003-486X},
   MRCLASS = {58F15 (58F11)},
  MRNUMBER = {1637655},
MRREVIEWER = {Nikolai Chernov},
       DOI = {10.2307/120960},
      URL = {https://doi.org/10.2307/120960},
}

@article {FFMA,
    AUTHOR = {Freitas, Ana Cristina Moreira and Freitas, Jorge Milhazes and
              Magalh\~{a}es, M\'{a}rio},
     TITLE = {Complete convergence and records for dynamically generated
              stochastic processes},
   JOURNAL = {Trans. Amer. Math. Soc.},
  FJOURNAL = {Transactions of the American Mathematical Society},
    VOLUME = {373},
      YEAR = {2020},
    NUMBER = {1},
     PAGES = {435--478},
      ISSN = {0002-9947},
   MRCLASS = {37A50 (37A25 37B20 60G55 60G70)},
  MRNUMBER = {4042881},
MRREVIEWER = {Johannes T. N. Krebs},
       DOI = {10.1090/tran/7922},
       URL = {https://doi.org/10.1090/tran/7922},
}

@article {Young2,
    AUTHOR = {Young, Lai-Sang},
     TITLE = {Recurrence times and rates of mixing},
   JOURNAL = {Israel J. Math.},
  FJOURNAL = {Israel Journal of Mathematics},
    VOLUME = {110},
      YEAR = {1999},
     PAGES = {153--188},
      ISSN = {0021-2172},
   MRCLASS = {37D25 (37A25 37D50)},
  MRNUMBER = {1750438},
MRREVIEWER = {Benoit Saussol},
       DOI = {10.1007/BF02808180},
       URL = {https://doi.org/10.1007/BF02808180},
}

@article {nonclt,
    AUTHOR = {Gou\"{e}zel, S\'{e}bastien},
     TITLE = {Central limit theorem and stable laws for intermittent maps},
   JOURNAL = {Probab. Theory Related Fields},
  FJOURNAL = {Probability Theory and Related Fields},
    VOLUME = {128},
      YEAR = {2004},
    NUMBER = {1},
     PAGES = {82--122},
      ISSN = {0178-8051},
   MRCLASS = {37A50 (37A30 37C30 37E05 60F05)},
  MRNUMBER = {2027296},
MRREVIEWER = {Jean-Ren\'{e} Chazottes},
       DOI = {10.1007/s00440-003-0300-4},
       URL = {https://doi.org/10.1007/s00440-003-0300-4},
}

@article {cuspsmorethanone,
    AUTHOR = {Jung, Paul and P\`ene, Francoise and Zhang, Hong-Kun},
     TITLE = {Convergence to {$\alpha$}-stable {L}\'evy motion for chaotic
              billiards with several cusps at flat points},
   JOURNAL = {Nonlinearity},
  FJOURNAL = {Nonlinearity},
    VOLUME = {33},
      YEAR = {2020},
    NUMBER = {2},
     PAGES = {807--839},
      ISSN = {0951-7715,1361-6544},
   MRCLASS = {37D40 (37A50 60F17)},
  MRNUMBER = {4048780},
MRREVIEWER = {Jianyu\ Chen},
       DOI = {10.1088/1361-6544/ab5148},
       URL = {https://doi.org/10.1088/1361-6544/ab5148},
}

\end{document}